\documentclass{smfart}
\usepackage[francais,english]{babel}
\usepackage{smfthm}
\usepackage{bull}
\theoremstyle{plain} 
\author{Pascal Redou}
\address{
Institut Girard Desargues, Universit\'e Claude Bernard Lyon 1, B\^atiment Braconnier (ex-101), 21 Avenue 
Claude Bernard, 69622 VILLEURBANNE Cedex, FRANCE 
}
\email{redou@enib.fr}
\title[Densit\'es tensorielles sur la sph\`ere]{Repr\'esentations de l'alg\`ebre de Lie conforme sur l'espace des densit\'es tensorielles sur la sph\`ere}
\alttitle{Representations of the conformal Lie algebra on the space of tensor densities on the sphere}

\newcommand{\bbR}{\mathbb{R}}

\newcommand{\bbC}{\mathbb{C}}

\newcommand{\Diff}{\mathrm{Diff}}
\newcommand{\Div}{\mathrm{Div}}

\newcommand{\End}{\mathrm{End}}
\newcommand{\cF}{{\mathcal{F}}}

\newcommand{\SL}{\mathrm{SL}}

\newcommand{\SO}{\mathrm{SO}}

\newcommand{\so}{\mathrm{o}}
\newcommand{\Tr}{\mathrm{Tr}}
\newcommand{\Vect}{\mathrm{Vect}}

\def\a{\alpha}

\def\d{\delta}

\def\l{\lambda}
\def\m{\mu}

\def\implies{\Rightarrow}

\def\mod#1{\left|{#1}\right|}

\begin{document}
\frontmatter

\begin{abstract}
\noindent Soit $\cF_\l(\mathbb{S}^n)$ l'espace des densit\'es tensorielles sur $\mathbb{S}^n$ de degr\'e $\l$
(ou, de fa\c{c}on \'equivalente, des densit\'es conformes de degr\'e $-\l n$). Cet espace est muni d'une structure de $\Diff(\mathbb{S}^n)$- et de 
$\Vect(\mathbb{S}^n)$-module, et 
nous d\'emontrons dans un premier temps qu'en tant que 
$\SO(n+1,1)$-module, il est infinit\'esimalement \'equivalent  au module 
 associ\'e
 \`a une repr\'esentation induite de la s\'erie principale du groupe connexe $\SO_0(n+1,1)$, plus pr\'ecis\'ement de la s\'erie sph\'erique nonunitaire sur l'espace ${\mathcal C}^\infty(\mathbb{S}^n)$. A partir des propri\'et\'es de cette derni\`ere, 
nous classifions selon les valeurs de $\l$ les sous-($\so(n+1,1),\SO(n+1))$-modules irr\'eductibles et unitaires de $\cF_\l(\mathbb{S}^n)$. 
\end{abstract}

\begin{altabstract}
\noindent Let $\cF_\l(\mathbb{S}^n)$ be the space of tensor densities on $\mathbb{S}^n$ with degree $\l$ (or, equivalently, of conformal densities 
with degree $\l$). This space is embedded with a structure of $\Diff(\mathbb{S}^n)$- and 
$\Vect(\mathbb{S}^n)$-module, and we first prove that as a $\SO(n+1,1)$-module, it is infinitesimally equivalent to an induced module from the 
principal series of the connected group 
$\SO_0(n+1,1)$, precisely from the nonunitary spherical series on the space ${\mathcal C}^\infty(\mathbb{S}^n)$. In function of $\l$, $(\mathfrak{g},K)$-simple and unitary submodules of 
$\cF_\l(\mathbb{S}^n)$ are therefore classified.
\end{altabstract}

\subjclass{G\'eom\'etrie diff\'erentielle et th\'eorie des repr\'esentations}
\keywords{Densit\'es tensorielles, $(\mathfrak{g}, K)$-modules, repr\'esentations induites.}
\altkeywords{Tensor densities, $(\mathfrak{g}, K)$-modules, induced representations}

\maketitle

\mainmatter
\section{Introduction et r\'esultats principaux}
L'espace $\cF_\l(M)$ des densit\'es tensorielles de degr\'e $\l$ sur une vari\'et\'e $M$ est l'espace des sections lisses du fibr\'e 
$\Delta_\l(M)=
\mod{\Lambda^nT^*M}^{\otimes\l}$ sur $M$. Cet espace est
utilis\'e dans les probl\`emes li\'es \`a la quantification g\'eom\'etrique, et plus r\'ecemment dans le domaine de la quantification
\'equivariante. Son introduction permet par exemple d'\'eviter l'\'ecueil de la non compatibilit\'e de l'action des op\'erateurs diff\'erentiels 
avec celle du groupe des diff\'eomorphismes de M (et donc de l'alg\`ebre de Lie $\Vect(M)$) sur l'espace ${\mathcal C}^\infty(M)$: l'exemple le plus 
c\'el\`ebre est l'espace des op\'erateurs de Sturm-Liouville $A=\frac{d^2}{dx^2}+u(x)$, qui agissent de $\cF_{-\frac{1}{2}}(\bbR)$ vers 
$\cF_{\frac{3}{2}}(\bbR)$. Il est de m\^eme aujourd'hui bien connu (voir par exemple \cite{ER}) que les seuls 
op\'erateurs diff\'erentiels lin\'eaires conform\'ement invariants (i.e. dont l'action commute avec celle de l'alg\`ebre de Lie conforme 
$\so(n+1,1)$) sur $\cF_\l(\bbR^n)$ sont les ``puissances du Laplacien'' :
pour tout $k\in\mathbb{N}$, il existe un unique (\`a constante pr\`es) op\'erateur diff\'erentiel lin\'eaire $\so(p+1,q+1)$-invariant d'ordre $k$ :
$$
A_{2k}: \cF_{\frac{n-2k}{2n}}
\to\cF_{\frac{n+2k}{2n}}
$$
et il n'existe pas d'autre op\'erateurs diff\'erentiels lin\'eaires $\so(p+1,q+1)$-invariants  d'ordre $k\geq1$ de $\cF_\l$ vers $\cF_\m$.

 De la m\^eme fa\c{c}on, la classification des op\'erateurs diff\'erentiels bilin\'eaires conform\'ement invariants a \'et\'e recemment effectu\'ee (\cite{OvR}), 
qui g\'en\'eralise la notion de transvectants ou de crochets de Rankin-Cohen.

\medskip

Notre travail s'inscrit \`a la crois\'ee de deux champs d'investigation : celui de la g\'eom\'etrie bien s\^ur, eu \'egard \`a la pr\'epond\'erance des densit\'es 
tensorielles dans le domaine de la g\'eom\'etrie diff\'erentielle, et, partant, de la physique th\'eorique, mais \'egalement celui de la th\'eorie des 
repr\'esentations, l'enjeu \'etant la classification exhaustive des espaces $\cF_\l(M)$ en tant que $\Diff(M)$-modules (sous-modules, modules \'equivalents, etc.) 

\medskip

Dans cet article, nous consid\'erons pr\'ecis\'ement  l'espace $\cF_\l(\mathbb{S}^n)$ des densit\'es tensorielles sur la sph\`ere $\mathbb{S}^n$. 
Le sens g\'eom\'etrique de cette \'etude r\'eside dans le fait que le groupe conforme $\SO(n+1,1)$ agit naturellement sur $\mathbb{S}^n$, ce qui conf\`ere
 \`a $\cF_\l(\mathbb{S}^n)$ une structure de 
$\SO(n+1,1)$-module, et par cons\'equent de $\so(n+1,1)$- et de $\SO(n+1)$-module. Nous ferons naturellement usage de la classification des 
repr\'esentations du groupe conforme,
 d\'ecrite par M. Takeshi Hirai (\cite{Hir}) ainsi que MM. Klimyk et Gavrilik (\cite{KG}).

\medskip

Donnons quelques pr\'ecisions sur les r\'esultats que nous allons d\'emontrer et les m\'ethodes employ\'ees : 

\medskip

 Rappelons que si $G$ est un groupe de Lie, $K$ un sous-groupe compact de $G$, et $\mathfrak{g}$ et $\mathfrak{k}$ sont leurs alg\`ebres de Lie 
respectives, on appelle 
$(\mathfrak{g},K)$-module un espace vectoriel complexe $E$ muni d'une repr\'esentation de $\mathfrak{g}$ et d'une repr\'esentation de $K$ v\'erifiant
\begin{enumerate}
\item $({\rm Ad}k\cdot X)\cdot e=k\cdot X\cdot k^{-1}\cdot e\quad \forall k\in K$, $X\in\mathfrak{g}$, $e\in E$
\item pour tout $e\in E$, $K\cdot e$ engendre un sous-espace vectoriel de dimension finie $F$ (on dit que $e$ est $K$-fini), la repr\'esentation de $K$ 
dans $F$ est continue et on a, pour tout $X\in\mathfrak{k}$,
$
X\cdot e=\displaystyle\frac{d}{dt}(\exp tX)\cdot e|_{t=0}$.
\end{enumerate}
Nous notons \`a pr\'esent $G=\SO_0(n+1,1)$, composante connexe de l'identit\'e dans $\SO(n+1,1)$, $\mathfrak{g}=\so(n+1,1)$, $K=\SO(n+1)$, 
et ${\mathcal H}(K)$ l'espace des vecteurs $K$-finis de $\cF_\l(\mathbb{S}^n)$.
Le r\'esultat final de notre travail est la classification, selon les valeurs de $\l$, des sous-$(\mathfrak{g},K)$-modules simples et unitaires de ${\mathcal H}(K)$.

\medskip

Dans un premier temps, nous d\'emontrons que la repr\'esentation du groupe conforme $\SO(n+1,1)$ sur $\cF_\l(\mathbb{S}^n)$ est 
infinit\'esimalement \'equivalente \`a une repr\'esentation induite de la s\'erie principale du groupe $G=\SO_0(n+1,1)$, plus pr\'ecis\'ement de la s\'erie dite 
sph\'erique nonunitaire :
notons $G=KAN$ la 
d\'ecomposition d'Iwasawa de $G$, $\rho$ la demi-somme des racines restreintes positives de la paire ($\so(n+1,1),\mathfrak{a})$, 
$A={\rm exp} \mathfrak{a}$.
Nous consid\'erons la repr\'esentation induite par le sous-groupe 
parabolique minimal $MAN$ de $G$, en faisant le choix de la repr\'esentation triviale du sous-groupe $M=\SO(n)$, centralisateur de $A$ 
dans $K$, et d'une repr\'esentation
 $\mu$ de dimension $1$ de $A$. Nous noterons, pour tout \'el\'ement $h$ de $A$,  $\mu(h)={\rm exp}(\nu(\log h))$, $\nu\in \mathfrak{a}^*$.  
Les notations assimilent 
un \'el\'ement de $\mathfrak{a}^*$ et sa valeur en l'\'el\'ement matriciel $H=E_{n+1,n+2}+E_{n+2,n+1}$. La d\'ecomposition d'Iwasawa montre que 
la repr\'esentation induite en question, que nous noterons 
${\rm Ind}_{MAN}^G(0\otimes \nu)$, prend effet sur l'espace des fonctions de ${\mathcal L}^2(K/M)={\mathcal L}^2(\mathbb{S}^n)$, et les op\'erateurs 
de cette repr\'esentation sont donn\'es, pour $g\in G$, par  
$$
{\rm Ind}_{MAN}^G(0\otimes \nu)(g)f(k)={\rm exp}(-\nu(\log h))f(k_g), {~\rm avec~} g^{-1}k=k_ghn\in KAN.
$$ 
Lorsque $\nu$ varie dans $\bbC$, nous obtenons ainsi les 
repr\'esentations de la s\'erie dite {\it sph\'erique nonunitaire}, qui munissent l'espace ${\mathcal L}^2(\mathbb{S}^n)$ d'une structure de $G$-module. 
Nous d\'esignons alors par ${\mathcal C}^\infty_\nu(\mathbb{S}^n)$ le sous-module form\'e des \'el\'ements de classe 
${\mathcal C}^\infty$. Nous d\'emontrons le r\'esultat-cl\'e suivant :

\begin{theo} 
$\cF_\l(\mathbb{S}^n)$ et ${\mathcal C}^\infty_\nu(\mathbb{S}^n)$ sont des $\mathfrak{g}$-modules isomorphes si et seulement si $\nu=n\l$, et cet 
isomorphisme est compatible avec l'action de $K$.
\label{equiv}
\end{theo} 
La d\'emonstration de ce r\'esultat pr\'esente l'int\'er\^et (et la difficult\'e!) d'exposer des calculs explicites sur un groupe de Lie, en l'occurence le groupe 
conforme, 
au moyen notamment de la d\'ecomposition d'Iwasawa. Elle s'appuie par ailleurs sur des consid\'erations cohomologiques simples 
mais d'une grande utilit\'e.
\medskip

En utilisant les propri\'et\'es des repr\'esentations de la s\'erie sph\'erique nonunitaire, nous classifions dans un second temps les 
sous-$(\mathfrak{g},K)$-modules simples et unitaires de ${\mathcal H}(K)$. Cette classification est la suivante :
\begin{theo}
\begin{enumerate}
\item Si $\l\not=l/n$ pour tout $l\in\mathbb{Z}$, ou si, pour $n>1$, $\l\in\{\frac{1}{n},\dots,\frac{n-1}{n}\}$, alors $\cF_\l(\mathbb{S}^n)$ 
contient un unique $(\mathfrak{g},K)$-module simple ${\mathcal H}(K)$, 
qui s'identifie 
\`a l'espace des polyn\^omes harmoniques  en les coordonn\'ees cart\'esiennes sur $\mathbb{S}^n$. Ce module est unitaire si et seulement 
si $\l=\frac{1}{2}+i\a$, $\a\in\bbR^*$, ou $\l\in ]0,1[\setminus\{\frac{1}{2}\}$.
\item Si $\l=-l/n$, $l\in\mathbb{N}$, ${\mathcal H}(K)$ admet un unique sous-$(\mathfrak{g},K)$-module simple, de dimension finie, form\'e de ses 
\'el\'ements 
de degr\'e $\le l$. Il est de plus unitaire pour $\l=0$.
\item Si $n=1$ et $\l=l$, $l\in\mathbb{N}^*$, ${\mathcal H}(K)$ admet deux sous-$(\mathfrak{g},K)$-modules simples, unitaires et de dimension infinie, 
dont la somme est l'ensemble des \'el\'ements de ${\mathcal H}(K)$ de degr\'e $\ge l$. 
\item Si $n>1$ et $\l=1+l/n$, $l\in\mathbb{N}$, ${\mathcal H}(K)$ admet un sous-$(\mathfrak{g},K)$-module simple, unitaire pour $\l=0$ et de dimension 
infinie, 
form\'e de ses \'el\'ements de degr\'e $\ge l+1$. 
\end{enumerate}
\label{class}
\end{theo}
\begin{rema}
Le r\'esultat obtenu contient \'egalement la classification des sous-quotients simples de l'espace des vecteurs $K$-finis, ce que la d\'emonstration mettra en exergue.

D'autre part, nous obtenons dans le m\^eme temps tous les sous-$G$-modules ferm\'es de $\cF_\l(\mathbb{S}^n)$ (d'apr\`es \cite{Kna}, th\'eor\`eme 8.9.), et 
la connexit\'e de $G$ implique que l'on obtient dans le cas (2) du th\'eor\`eme \ref{class} tous les sous $\mathfrak{g}$-modules simples de dimension finie de 
$\cF_\l(\mathbb{S}^n)$.
\end{rema}

Nous commencerons par pr\'esenter l'espace des densit\'es tensorielles sur la sph\`ere, puis la s\'erie sph\'erique nonunitaire du groupe 
$\SO_0(n+1,1)$, pour entreprendre ensuite la d\'emonstration des th\'eor\`emes \ref{equiv} et \ref{class}.
\section{Densit\'es tensorielles}
\subsection{Structures de $\mathfrak{g}$- et $K$-module}
Soit $\cF_\l(\mathbb{S}^n)$ l'espace des densit\'es tensorielles de degr\'e $\l\in\bbC$ sur la sph\`ere $\mathbb{S}^n$. Cet espace est de fa\c{c}on 
naturelle muni d'une structure de $\Diff(\mathbb{S}^n)$- et $\Vect(\mathbb{S}^n)$-module. En tant qu'espace vectoriel, 
$\cF_\l(\mathbb{S}^n)$ est isomorphe \`a l'espace de fonctions lisses \`a valeurs complexes ${\mathcal C}^\infty_\bbC(\mathbb{S}^n)$, mais l'action d'un champ de vecteurs 
$X=\sum_{i=1}^n X^i\frac{\partial}{\partial{}x_i}$ de l'alg\`ebre de Lie
 $\Vect(\mathbb{S}^n)$ est donn\'ee par la d\'eriv\'ee de Lie de degr\'e $\l$
\begin{equation}
L^\l_{X}=
X^i\frac{\partial}{\partial x_i} +\l\,\Div(X),
\label{LieDer}
\end{equation}
ind\'ependamment du choix du syst\`eme de coordonn\'ees.
Explicitement, si $\varphi\in\cF_\l(\mathbb{S}^n)$, et si nous notons $\Omega$ la forme de volume 
$$
\Omega=dx_1\wedge\cdots\wedge dx_n,
$$
la formule (\ref{LieDer}) se lit
$$
L_X^\l(\varphi(x_1,\ldots,x_n)(\Omega)^\l)=\sum_{i=1}^n(X^i\frac{\partial{} \varphi}{\partial{}x_i}+
\l\frac{\partial {}X^i}{\partial{}x_i}\varphi)(x_1,\ldots,x_n)(\Omega)^\l~.
$$

\noindent L'alg\`ebre de Lie $\so(n+1,1)$ des transformations conformes infinit\'esimales, que nous d\'esignons plus bri\`evement par alg\`ebre de 
Lie conforme, 
est engendr\'ee par les champs de vecteurs
\begin{equation}
\begin{array}{rcl}
X_i  
&=&
\displaystyle \frac{\partial}{\partial s_i}\;\hfill\cr
\noalign{\smallskip}
X_{ij}
&=&
\displaystyle s_i\frac{\partial}{\partial s_j}-
s_j\frac{\partial}{\partial s_i}\;\hfill\cr
\noalign{\smallskip}
X_0  
&=&
\sum_i\displaystyle s_i\frac{\partial}{\partial s_i}\;\hfill\cr
\noalign{\smallskip}
\bar X_i 
&=&
\sum_j(\displaystyle s_j^2\frac{\partial}{\partial s_i}-
2s_is_j\frac{\partial}{\partial s_j})\hfill\cr
\end{array}
\label{confGenerators}
\end{equation}
o\`u $(s_1,\ldots,s_n)$ sont les coordonn\'ees {\it st\'er\'eographiques} sur la sph\`ere $\mathbb{S}^n$. 

\medskip

\noindent Nous consid\'erons l'espace $\cF_\l(\mathbb{S}^n)$ comme un $\mathfrak{g}=\so(n+1,1)$-module, au moyen de la formule (\ref{LieDer}) et 
des champs de vecteurs 
(\ref{confGenerators}).
D'autre part, nous d\'eduisons de l'action du groupe $\Diff(\mathbb{S}^n)$ celle du sous-groupe $K=\SO(n+1)$, donn\'ee 
par la repr\'esentation dite r\'eguli\`ere \`a gauche
\begin{equation}
(k_0\cdot f)(k)=f(k_0^{-1}k),
\label{reg}
\end{equation} 
o\`u $k_0\in K$, $k\in\mathbb{S}^n\simeq\SO(n+1)/\SO(n)$.
 La structure de $K$-module de $\cF_\l(\mathbb{S}^n)$ est donc 
ainsi explicit\'ee.

\medskip

 Nous d\'ecrivons ci-apr\`es comment nous pouvons utiliser les coordonn\'ees cart\'esiennes
pour d\'ecrire la structure de $\mathfrak{g}$-module de $\cF_\l(\mathbb{S}^n)$, ce dont nous ferons usage par la suite.
\subsection{Module des fonctions homog\`enes sur $\bbR^{n+1}\setminus\{0\}$}
Nous notons ${\it C}^\infty_{-\l(n+1)}(\bbR^{n+1}\setminus\{0\})$ l'espace des fonctions lisses sur
 $\bbR^{n+1}\setminus\{0\}$, homog\`enes de degr\'e $-(n+1)\l$, i.e. telles que
$$
f(\alpha x_0,\ldots,\alpha x_n)~=~\alpha^{-(n+1)\l}f(x_0,\ldots,x_n),~ \forall\alpha\in\bbC.
$$
Cet espace est naturellement muni d'une structure de ${\rm Vect}(\bbR^{n+1})$-module, {\it l'action d'un vecteur \'etant donn\'ee par la d\'eriv\'ee de 
Lie le long de ce vecteur}. 
D\'esignons par ${\it C}^\infty_{-\l(n+1)}(\bbR^{n+1}\setminus\{0\}|_{\mathbb{S}^n})$ la restriction de cet espace aux fonctions sur $\mathbb{S}^n$.
Nous lui attribuons une structure de $\Diff(\mathbb{S}^n)$-module, h\'erit\'ee de celle de $\cF_\l(\mathbb{S}^n)$, de la mani\`ere suivante :

\begin{prop} (\cite{Ovs})
L'espace $\cF_\l(\mathbb{S}^n)$ est isomorphe en tant que $\Diff(\mathbb{S}^n)$-module \`a l'espace 
${\it C}^\infty_{-\l(n+1)}(\bbR^{n+1}\setminus\{0\}|_{\mathbb{S}^n})$.
\label{Reltheo}
\end{prop}
\noindent Explicitement, si $\varphi\in\cF_\l(\mathbb{S}^n)$ et si nous notons $t_1,\ldots t_n$ les coordonn\'ees {\it projectives} sur $\mathbb{S}^n$ et 
$x_0,x_1,\dots ,x_n$ les coordonn\'ees {\it cart\'esiennes} sur $\bbR^{n+1}$, cet isomorphisme, que nous appelons {\it rel\`evement}, est donn\'e par
\begin{equation}
\displaystyle\varphi(t_1,\dots t_n)(dt_1\wedge\cdots\wedge dt_n)^{\l}\longmapsto\varphi(\frac{x_1}{x_0},\dots ,\frac{x_n}{x_0})(x_0)^{-(n+1)\l}.
\label{rel}
\end{equation}

Bien entendu, il correspond \`a ce rel\`evement un ``rel\`evement des champs de vecteurs", faisant correspondre aux champs (\ref{confGenerators}) 
qui agissent par la d\'eriv\'ee de Lie de degr\'e $\l$ sur 
$\cF_\l(\mathbb{S}^n)$, des champs de vecteurs en coordonn\'ees cart\'esiennes, qui agissent par la d\'eriv\'ee de Lie standard (degr\'e $0$) sur 
${\it C}^\infty_{-\l(n+1)}(\bbR^{n+1}\setminus\{0\})$. N'en ayant pas ici l'utilit\'e, nous n'explicitons pas ces~champs.
\subsection{Polyn\^omes harmoniques}
Par la suite, il nous sera n\'ecessaire de nous placer alternativement dans chacun des espaces 
${\it C}^\infty_{-\l(n+1)}(\bbR^{n+1}\setminus\{0\}|_{\mathbb{S}^n})$ et $\cF_\l(\mathbb{S}^n)$, utilisant \`a cette fin l'isomorphisme que nous venons de d\'ecrire. Nous ferons en 
particulier un large 
usage des polyn\^omes harmoniques dans $\bbR^{n+1}|_{\mathbb{S}^n}$, i.e des polyn\^omes en les $n+1$ variables r\'eelles
$x_0,\ldots,x_n$, restreints \`a la sph\`ere et appartenant au noyau du Laplacien standard       
$$
\Delta=\displaystyle\frac{\partial^2}{\partial x_0^2}+\cdots+\frac{\partial^2}{\partial x_n^2}.
$$
Parmi ces polyn\^omes, nous consid\'erons en particulier les \'el\'ements 
$$
\psi_{il}~=~(x_0+{\bf i}x_i)^l,
$$
o\`u nous avons not\'e ${\bf i}=\sqrt{-1}$. Afin d'expliciter un exemple concret, mais aussi parce que nous en aurons besoin, montrons comment 
$\psi_{il}$ peut \^etre consid\'er\'e comme un
 \'el\'ement       
de l'espace  $\cF_\l(\mathbb{S}^n)$: \`a l'aide de la relation $x_0^2+\cdots+x_n^2=1$, nous obtenons
$$
\psi_{il}~=~\displaystyle\frac{(x_0+{\bf i}x_i)^l}{(x_0^2+\cdots+x_n^2)^{\frac{l}{2}+\l\frac{n+1}{2}}},
$$
et v\'erifions facilement que $\psi_{il}\in{\mathcal C}^\infty_{-(n+1)\l}(\bbR^{n+1}\setminus\{0\}|_{\mathbb{S}^n})$.
\begin{rema}
L'action d'un \'el\'ement de l'alg\`ebre de Lie $\so(n+1)$ sur $\psi_{il}$ ne fait pas intervenir $\lambda$: il suffit d'appliquer un op\'erateur 
infinit\'esimal 
$x_k\displaystyle\frac{\partial }{\partial x_l}-x_l\displaystyle\frac{\partial }{\partial x_k}$ de $\so(n+1)$ pour s'en convaincre. 
\end{rema}
\noindent Utilisant maintenant l'isomorphisme (\ref{rel}), nous obtenons l'image de $\psi_{il}$ dans $\cF_\l(\mathbb{S}^n)$:
$$
\psi_{il}~=~\displaystyle\frac{(1+{\bf i}t_i)^l}{(1+t_1^2+\cdots+t_n^2)^{\frac{l}{2}+\l\frac{n+1}{2}}}(dt_1\wedge\cdots\wedge dt_n)^\l 
\label{psiil}
$$
(par commodit\'e nous la d\'esignons aussi par $\psi_{il}$), cette expression \'etant \'evidemment donn\'ee en coordonn\'ees projectives. 
Il nous reste alors \`a appliquer les champs (\ref{confGenerators}) exprim\'es dans ces coordonn\'ees, si nous souhaitons obtenir les images de 
$\psi_{il}$
 par l'action des champs de $\so(n+1,1)$. Utilisant \`a nouveau (\ref{rel}), nous obtiendrons les expressions de ces images dans l'espace 
${\it C}^\infty_{-\l(n+1)}(\bbR^{n+1}\setminus\{0\}|_{\mathbb{S}^n})$.

\section{Repr\'esentations de la s\'erie sph\'erique nonunitaire}
L'id\'ee principale de notre travail de classification r\'eside dans l'identification de $\cF_\l(\mathbb{S}^n)$ comme module induit : il s'agit d'un 
module de la s\'erie principale des repr\'esentations du groupe $\SO_0(n+1,1)$, plus pr\'ecis\'ement de la s\'erie dite sph\'erique nonunitaire. 
Nous d\'ecrivons ci-apr\`es la construction de cette s\'erie (que l'on trouvera expos\'ee dans l'ouvrage  \cite{ViK} de MM.~Klymik et Vilenkin). 
\subsection{Description}
Avant de consid\'erer plus pr\'ecis\'ement le groupe conforme, nous donnons quelques r\'esultats g\'en\'eraux.
Soit $G$ un groupe de Lie lin\'eaire non compact, $K$ son sous-groupe compact maximal. On note $\mathfrak{g}$ l'alg\`ebre de Lie de $G$, et 
$\mathfrak{k}$ celle de $K$, qui est une sous-alg\`ebre de $\mathfrak{g}$.
Il existe dans $\mathfrak{g}$ un isomorphisme involutif $\theta$, l'involution de Cartan, pour lequel $\mathfrak{k}$ est le sous-espace stationnaire. 
Le sous-espace 
$\{X:\theta X=-X\}$ dans $\mathfrak{g}$ est not\'e $\mathfrak{p}$. Nous avons alors 
$$
\mathfrak{g}=\mathfrak{k}+\mathfrak{p}.
$$
Cette d\'ecomposition se transforme par l'application exponentielle $\mathfrak{g}\longrightarrow G$, en la d\'ecomposition $G=KP$ de $G$, o\`u 
$P=\exp\mathfrak{p}$.

Soit $\mathfrak{a}$ une sous-alg\`ebre commutative maximale de $\mathfrak{p}$. La dimension de $\mathfrak{a}$ est appel\'ee le rang de 
$\mathfrak{g}$ 
et de $G$. Le sous-groupe $$A=\exp\mathfrak{a}$$ est commutatif.

Les op\'erateurs de repr\'esentation adjointe ${\rm ad}H$, $H\in\mathfrak{a}$, sont anti-hermitiens pour le produit 
$$
\langle X,Y\rangle=-\Tr({\rm ad}X{\rm ad}(\theta Y)),~X,Y\in\mathfrak{g},
$$
et par cons\'equent nous avons la d\'ecomposition
\begin{equation}
\mathfrak{g}=\mathfrak{g}_0+\sum_\gamma\mathfrak{g}_\gamma~,
\label{dec1}
\end{equation}
o\`u $\mathfrak{g}_0$ est le noyau de l'op\'erateur ${\rm ad}H$ et $\mathfrak{g}_\gamma$ correspond aux valeurs propres $\gamma(H)$, 
$H\in\mathfrak{a}$. La d\'ecomposition (\ref{dec1}) est orthogonale. Les formes lin\'eaires  $\gamma$ sont appel\'ees {\it racines restreintes} 
de la paire $(\mathfrak{g},\mathfrak{a})$,
et les sous-espaces $\mathfrak{g}_\gamma$ sont les {\it sous-espaces radiciels}.

Si $H_1,\ldots H_l$ est une base de $\mathfrak{a}$, et  si le premier nombre non nul dans la suite
\newline
$\{\gamma(H_1),\ldots,\gamma(H_l)\}$ est positif (resp. n\'egatif),
alors la racine $\gamma$ est dite positive (resp. n\'egative). La dimension de $\mathfrak{g}_\gamma$ est appel\'ee multiplicit\'e de $\gamma$, 
et est not\'ee 
$m(\gamma)$.
Posons
\begin{equation}
\rho=\displaystyle\frac{1}{2}\sum_{\gamma>0}m(\gamma)\gamma.
\label{ro}
\end{equation}
La somme 
$$
\mathfrak{n}=\sum_{\gamma>0}\mathfrak{g}_\gamma
$$
est une sous-alg\`ebre nilpotente maximale de $\mathfrak{g}$, et
$$
N=\exp\mathfrak{n}
$$
est un sous-groupe nilpotent maximal de $G$.

\noindent Le groupe $G$ admet une d\'ecomposition d' {\it Iwasawa} 
$$
G=KAN,
$$
ce qui signifie que chaque \'el\'ement $g\in G$ peut \^etre \'ecrit de mani\`ere unique $g=khn$, o\`u $k\in K$, $h\in A,~n\in N$.
De plus l'application $(k,h,n)\mapsto khn$ est un diff\'eomorphisme analytique de $K\times A\times N$ vers $G$.

Soit $M$ le centralisateur de $A$ dans $K$. Le sous-groupe 
$$
P=MAN
$$
est appel\'e sous-groupe parabolique minimal de $G$. Un sous-groupe $P'$ contenant $P$ et diff\'erent de $G$ est appel\'e sous-groupe parabolique. 
Un tel groupe s'obtient par extension de $P=MAN$ en $P'=M'AN$, o\`u $M\subset M'\subset K$.

\bigskip

 Soient $G,K,A,M,N,P$ les groupes ainsi d\'efinis. Nous choisissons une repr\'esentation de dimension $1$
$$
\mu(h)=\exp(\nu(H)),\quad h\in \exp H,
$$
du sous-groupe $A=\exp \mathfrak{a}$. Alors la correspondance
$$
p=mhn\longmapsto\d(mhn)=\mu(h)
$$
d\'efinit  une repr\'esentation de dimension $1$ du sous-groupe parabolique minimal $P$.
Cette repr\'esentation induit une repr\'esentation ${\rm Ind}_{MAN}^G(0\otimes \nu)$ du groupe $G$, que nous noterons par souci de concision 
$I_\nu$, agissant 
sur l'espace des fonctions $f$ sur $G$ qui satisfont \`a la condition
\begin{equation}
f(gp)=\mu(h^{-1})f(g),\quad p=mhn \in P.
\label{ind1}
\end{equation}
Les op\'erateurs $I_\nu(g_0)$, $g_0\in G$, agissent sur ces fonctions par la formule 
$$
I_\nu(g_0)f(g)=f(g_0^{-1}g).
$$
Nous appelons ces repr\'esentations ${\rm Ind}_{MAN}^G(0\otimes \nu)=I_\nu$ {\it repr\'esentations de la s\'erie sph\'erique nonunitaire}.

Si des fonctions $f$ sur $G$ satisfont \`a la condition (\ref{ind1}), elles sont d\'etermin\'ees par leurs valeurs sur un certain sous-groupe de $G$ :
 la d\'ecomposition d'Iwasawa $G=KAN$ montre qu'elles sont d\'etermin\'ees par leurs valeurs sur $K$, et la relation
\begin{equation}
f(km)=f(k),\quad m\in M,
\label{quot}
\end{equation}
est satisfaite.

Les op\'erateurs $I_\nu(g)$ sont donn\'es par la formule 
$$
I_\nu(g)f(k)=\mu(h^{-1})f(k_g),
$$
o\`u $h\in A$ et $k_g\in K$ sont d\'efinis par la d\'ecomposition d'Iwasawa
$$
g^{-1}k=k_ghn
$$
de l'\'el\'ement $g^{-1}k$.

\medskip

\noindent Le produit scalaire
$$
\langle f_1,f_2\rangle=\int_K f_1(k)\overline{f_2(k)}dk
$$
est introduit dans l'espace des fonctions sur $K$.
La relation (\ref{quot}) signifie que ces fonctions sont en fait des fonctions sur le groupe quotient $K/M$.

Les repr\'esentations $I_\nu$ sont par cons\'equent r\'ealis\'ees dans l'espace des fonctions sur $K/M$.
\subsection{D\'ecomposition du groupe $\SO_0(n+1,1)$}
Pour ce groupe, nous avons $K\cong\SO(n+1)$: tout \'el\'ement de $K$ s'\'ecrit en effet sous la forme
$$
k=
\left[
\begin{array}{c|c}
 \begin{array}{ccc}
\hspace{0.2cm}&\hspace{0.2cm}&\hspace{0.2cm}\\
&R&\\
&&
\end{array}&~{\bf 0}\\
\hline
{\bf 0}& ~1
\end{array}
\right]
$$
o\`u $R\in\SO(n+1)$. 

\noindent D'autre part,  le sous-groupe $A$ est constitu\'e des matrices
\begin{equation}
h_{n+1}(t)=
\left[
\begin{array}{c|c}
 \begin{array}{ccc}
\hspace{0.3cm}&\hspace{0.3cm}&\hspace{0.3cm}\\
&I_n&\\
&&
\end{array}&~{\bf 0}\\
\hline
{\bf 0}& \begin{array}{cc}
	\cosh t&\sinh t\\
	\sinh t&\cosh t
	\end{array}
\end{array}
\right]
,
\label{A}
\end{equation}
o\`u $I_n$ est la matrice identit\'e de type $(n,n)$.

Nous noterons $H$ l'\'el\'ement de $\so(n+1,1)$ tel que $h_{n+1}(t)={\rm exp}tH$.
Ainsi, 
\begin{equation}
H=E_{n+1,n+2}+E_{n+2,n+1}, 
\label{h}
\end{equation}
o\`u $E_{ij}$ est la matrice dont les coefficients sont les $(E_{ij})_{rs}=\d_{ir}\d_{js}$.

Le groupe $M$ est isomorphe \`a $\SO(n)$. Le groupe $N$ est constitu\'e des matrices 
\begin{equation}
n({\bf a})=
\left[
\begin{array}{c|cc}
 \begin{array}{ccc}
\hspace{0.5cm}&\hspace{0.5cm}&\hspace{0.5cm}\\
&I_n&\\
&&
\end{array}
&
{\bf a}^t&-{\bf a}^t
\\
\hline
-{\bf a}&1-\frac{\|{\bf a}\|^2}{2}&\frac{\|{\bf a}\|^2}{2}\\
-{\bf a}&-\frac{\|{\bf a}\|^2}{2}&1+\frac{\|{\bf a}\|^2}{2}
\end{array}
\right]
,
\label{N}
\end{equation}
o\`u ${\bf a}=(a_1,\ldots,a_n)\in\bbR^n$, et $||{\bf a}||^2=\sum a_j^2$.
Nous v\'erifions facilement que $N$ est ici ab\'elien.

\medskip

En dernier lieu, pour $g\in\SO_0(n+1,1)$, le param\`etre $t$ de $h_{n+1}(t)\in A$, dans la d\'ecomposition d'Iwasawa 
$g=kh_{n+1}(t)n$, $g\in\SO_0(n+1,1)$, 
est donn\'e par $t=\ln[{\bf x_0},g{\bf \xi_0}]$, o\`u $
[{\bf x},{\bf y}]=x^1y^1+\cdots+x^{n+1}y^{n+1}-x^{n+2}y^{n+2}$, ${\bf \xi_0}=(0,\ldots ,0,1,1)$, et ${\bf x_0}=(0,\ldots,0,1)$.
\subsection{S\'erie sph\'erique nonunitaire du groupe $\SO_0(n+1,1)$}
Nous avons ici
$$
\mathfrak{a}=\{t(E_{n+1,n+2}+E_{n+2,n+1}),~t\in\bbR\},
$$
ainsi les caract\`eres du sous-groupe $A$ et, partant, les repr\'esentations ${\rm Ind}_{MAN}^G(0\otimes \nu)=I_\nu$ de la s\'erie sph\'erique nonunitaire 
du groupe $\SO_0(n+1,1)$, 
sont-ils donn\'es par un nombre complexe $\nu=\nu(H)$, o\`u $H$ est donn\'e par (\ref{h}).
On notera donc ult\'erieurement $\nu(H)=\nu$ et $\rho(H)=\rho$, $\rho$ \'etant d\'efini par (\ref{ro}). Nous avons alors en particulier $\rho=\frac{n}{2}$.

Puisque, pour $\SO_0(n+1,1)$, nous avons $K=\SO(n+1)$ et $M=\SO(n)$, les repr\'esentations $I_\nu$ sont r\'ealis\'ees dans l'espace
 ${\mathcal L}^2(\mathbb{S}^n)$ des fonctions de carr\'e int\'egrable sur la sph\`ere 
$$
\mathbb{S}^n=\SO(n+1)/\SO(n)\subset\bbR^{n+1}.
$$
De telles fonctions peuvent \^etre consid\'er\'ees comme fonctions des coordonn\'ees sph\'eriques $\theta_1,\dots,\theta_n$ sur $\mathbb{S}^n$. 
En effet consid\'erons l'\'el\'ement $P$ de $\mathbb{S}^n$, de coordonn\'ees sph\'eriques $(\theta_1,\dots,\theta_n)$.
 Nous avons alors en coordonn\'ees cart\'esiennes $(x_0,\dots,x_n)$
$$
P\left(
	\begin{array}{c}
	\sin\theta_1\cdots\sin\theta_n\\
	\cos\theta_1\sin\theta_2\cdots\sin\theta_n\\	 
	\cos\theta_2\sin\theta_3\cdots\sin\theta_n\\
	\vdots\\
	\cos\theta_{n-1}\sin\theta_n\\
	\cos\theta_n
	\end{array}
  \right).
$$
De plus, l'application
\begin{equation}
\varphi :\SO(n+1)\longrightarrow\bbR^{n+1},\quad\varphi(g)=ge_{n+1},
\label{isosphere}
\end{equation}
o\`u $e_{n+1}$ a pour coordonn\'ees cart\'esiennes $(0,\ldots,0,1)$,
donne le diff\'eomorphisme 
$$
\SO(n+1)/\SO(n)\cong\mathbb{S}^n,
$$
puisque le stabilisateur de $e_{n+1}$ dans  $\SO(n+1)$ est donn\'e par les matrices 
$$
\left[
\begin{array}{c|c}
 \begin{array}{ccc}
&&\\
&R\in\SO(n)&\\
&&
\end{array}&~{\bf 0}\\
\hline
{\bf 0}& ~1
\end{array}
\right]
.
$$
Dans cet isomorphisme, \`a l'\'el\'ement $P\in \mathbb{S}^n$ donn\'e ci-dessus correspond la classe
$$
\left[
\begin{array}{c|c}
 \begin{array}{ccc}
\hspace{0.2cm}&\hspace{0.2cm}&\hspace{0.2cm}\\
&R\in\SO(n)&\\
&&
\end{array}
&
\begin{array}{c}\sin\theta_1\cdots\sin\theta_n\\
\cos\theta_1\sin\theta_2\cdots\sin\theta_n\\	 
	\cos\theta_2\sin\theta_3\cdots\sin\theta_n\\
 \vdots\\ \cos\theta_{n-1}\sin\theta_n\end{array}\\
\hline

{\bf 0}&\cos\theta_n
\end{array}
\right]
\in \SO(n+1)/\SO(n)
.
$$

\bigskip

Nous allons maintenant d\'eterminer l'expression 
\begin{equation}
I_\nu(g)f(k)=\mu(h^{-1})f(k_g),
\label{ind4}
\end{equation}
o\`u $k\in K/M\cong \mathbb{S}^n,\quad g^{-1}k=k_ghn \in KAN,$
pour $g\in A$ et $g\in N$.

\medskip

Notons $k_i$ la rotation d'angle $\theta_i$ dans le plan $(e_i,e_{i+1})$, donn\'ee par 
$$
k_i=
\left[
\begin{array}{cccccc}
1&0&&&&0\\
0&\ddots&&&&0\\
&&\cos\theta_i&\sin\theta_i &&\\
&&-\sin\theta_i&\cos\theta_i&&\\
&&&\ddots&\\
0&&&&0&1\\
\end{array}
\right]
\in K=\SO(n+1)
.
$$
De m\^eme, $k'_i$ d\'esignera une semblable rotation, d'angle $\theta'_i$.

\noindent A l'\'el\'ement $P\in\mathbb{S}^n$ pr\'esent\'e plus haut correspond la classe d'\'equivalence dans $\SO(n+1)/SO(n)$ de l'\'el\'ement 
$k_1k_2\cdots k_n$, puisque la $n+1$-\`eme colonne de l'\'el\'ement $k_1k_2\cdots k_n$ est pr\'ecis\'ement form\'ee des coordonn\'ees cart\'esiennes 
de $P$.
Nous avons en effet (comme \'el\'ement du groupe $\SO_0(n+1,1)$),
$$
\begin{array}{l}
k_1k_2\cdots k_n=\\[15pt]
\left[
\begin{array}{cccccc}
\cos\theta_1&\sin\theta_1\cos\theta_2&\ldots&\sin\theta_1\cdots\sin\theta_{n-1}\cos\theta_n&\sin\theta_1\cdots\sin\theta_n&0\\
-\sin\theta_1&\cos\theta_1\cos\theta_2&\ldots&\cos\theta_1\sin\theta_2\cdots\sin\theta_{n-1}\cos\theta_n&\cos\theta_1\sin\theta_2\cdots\sin\theta_n&0\\
0&-\sin\theta_2&\ldots&\cos\theta_2\sin\theta_3\cdots\sin\theta_{n-1}\cos\theta_n&\cos\theta_2\sin\theta_3\cdots\sin\theta_n&0\\
\vdots&\vdots&&\vdots&\vdots&\vdots\\
0&0&\ldots&-\sin\theta_n&\cos\theta_n&0\\
0&&&&0&1
\end{array}
\right].
\end{array}
$$

\bigskip

$\bullet$ Consid\'erons tout d'abord $g\in A$, $g=h_{n+1}(t)$ comme dans ($\ref{A}$), et posons $k=k_1k_2\cdots k_n\in K/M$.
Utilisant alors la d\'ecomposition d'Iwasawa et le fait que $g\in A$ commute avec chaque \'el\'ement $k_ i$ appartenant \`a $M=\SO(n)$, 
nous obtenons
$$
g^{-1}k=k_ghn,\quad h=h_{n+1}(u),~n\in N,
$$
o\`u $k_g$ est du type
$$
k_g=k_1k_2\cdots k_{n-1}k'_n,
$$
$k'_n$ \'etant caract\'eris\'ee par l'angle $\theta'_n$, tel que 
\begin{equation}
\cos \theta'_n=\displaystyle\frac{\cos \theta_n\cosh t-\sinh t}{\cosh t-\cos \theta_n\sinh t}.
\label{teta'}
\end{equation}
D'autre part une identification matricielle directe nous donne $g^{-1}k=k_ghn$, $h=h_{n+1}(u)$, o\`u
$$
e^u=\cosh t-\cos \theta_n\sinh t.
$$
Enfin, au moyen de l'isomorphisme donn\'e par (\ref{isosphere}), nous obtenons la relation
\begin{equation}
I_\nu(h_{n+1}(t))f(\theta_1,\dots,\theta_n)=(\cosh t-\cos \theta_n\sinh t)^{-\nu} f(\theta_1,\dots,\theta_{n-1},\theta'_n),
\label{indA}
\end{equation}
o\`u $\theta'_n$ est donn\'e par (\ref{teta'}).

\bigskip

$\bullet$ Consid\'erons \`a pr\'esent le cas $g\in N$, $g=n({\bf a})$, comme dans (\ref{N}) :

\noindent \`A nouveau nous posons $k=k_1k_2\cdots k_n\in K/M$, et nous \'ecrivons la d\'ecomposition d'Iwasawa pour $g^{-1}k$. Chaque \'el\'ement  
du groupe 
quotient $K/M=\SO(n+1)/\SO(n)$ est du type 
$$
k_g=k'_1k'_2\cdots k'_nm,\quad m\in \SO(n).
$$
Un calcul matriciel direct nous conduit \`a l'identification
$$
g^{-1}k=k_ghn\iff n({\bf a})^{-1}k_1k_2\cdots k_n=k'_1k'_2\cdots k'_nmh_{n+1}(u)n({\bf \alpha}),
$$
o\`u ${\bf\alpha}=(\alpha_1,\dots,\alpha_n)$. Il convient de remarquer que $u$ et $\theta'_i$, $1\leq i\leq n$ doivent \^etre d\'etermin\'es, 
afin d'obtenir l'expression (\ref{ind4}).
L'identification ci-dessus nous conduit aux relations:
$$
\begin{array}{rcl}
e^u&=&1+\displaystyle\frac{||{\bf a}||^2}{2}+a_1\sin\theta_1\cdots\sin\theta_n + a_2\cos\theta_1\sin\theta_2\cdots\sin\theta_n \\[6pt]
&&+\cdots+
a_n\cos\theta_{n-1}\sin\theta_n-
\displaystyle\frac{||{\bf a}||^2}{2}\cos\theta_n;\\[15pt]
\cos\theta'_n&=&e^{-u}\Big(\displaystyle\frac{||{\bf a}||^2}{2}+a_1\sin\theta_1\cdots\sin\theta_n + a_2\cos\theta_1\sin\theta_2\cdots\sin\theta_n\\[6pt]&&
\hspace{1cm}+\cdots+
a_n\cos\theta_{n-1}\sin\theta_n
+(1-\displaystyle\frac{||{\bf a}||^2}{2})\cos\theta_n\Big);\\[15pt]
&=&1+\displaystyle\frac{\cos\theta_n -1}{e^u};
\end{array}
$$
$$
\begin{array}{rcl}
\cos\theta'_{n-1}\sin\theta'_n&=&e^{-u}(a_n+\cos\theta_{n-1}\sin\theta_n-a_n\cos\theta_n);\\[15pt]
\cos\theta'_{n-2}\sin\theta'_{n-1}\sin\theta'_n&=&e^{-u}(a_{n-1}+\cos\theta_{n-2}\sin\theta_{n-1}\sin\theta_n-a_{n-1}\cos\theta_n);\\
&\vdots&\\
\cos\theta'_1\sin\theta'_2\cdots\sin\theta'_n&=&e^{-u}(a_2+\cos\theta_1\sin\theta_2\cdots\sin\theta_n-a_2\cos\theta_n)\\[15pt]
\sin\theta'_1\sin\theta'_2\cdots\sin\theta'_n&=&e^{-u}(a_1+\sin\theta_1\sin\theta_2\cdots\sin\theta_n-a_1\cos\theta_n).
\end{array}
$$
Ces \'egalit\'es nous m\`enent enfin aux relations de r\'ecurrence suivantes:
\begin{equation}
\begin{array}{rcl}
\tan\theta'_1&=&\displaystyle\frac{a_1+\sin\theta_1\sin\theta_2\cdots\sin\theta_n-a_1\cos\theta_n}{a_2+
\cos\theta_1\sin\theta_2\cdots\sin\theta_n-a_2\cos\theta_n};\\[15pt]
\cos\theta'_1\tan\theta'_2&=&\displaystyle\frac{a_2+\cos\theta_1\sin\theta_2\cdots\sin\theta_n-a_2\cos\theta_n}{a_3+\cos\theta_2\sin\theta_3\cdots\sin\theta_n
-a_3\cos\theta_n};\\
&\vdots&\\
\cos\theta'_{n-2}\tan\theta'_{n-1}&=&\displaystyle\frac{a_{n-1}+\cos\theta_{n-2}\sin\theta_{n-1}\sin\theta_n-a_{n-1}\cos\theta_n}
{a_n+\cos\theta_{n-1}\sin\theta_n
-a_n\cos\theta_n};\\[15pt]
\cos\theta'_n&=&1+\displaystyle\frac{\cos\theta_n -1}{e^u};
\label{thet'exp1}
\end{array}
\end{equation}
\begin{equation}
\begin{array}{rcl}
e^u&=&1+\displaystyle\frac{||{\bf a}||^2}{2}+a_1\sin\theta_1\cdots\sin\theta_n + a_2\cos\theta_1\sin\theta_2\cdots\sin\theta_n \\[6pt]
&&+\cdots+
a_n\cos\theta_{n-1}\sin\theta_n-
\displaystyle\frac{||{\bf a}||^2}{2}\cos\theta_n.
\label{thet'exp2}
\end{array}
\end{equation}

\noindent En cons\'equence, l'action d'un \'el\'ement du groupe nilpotent $N$ est donn\'ee par

\begin{equation}
I_\nu(n({\bf a}))f(\theta_1,\dots,\theta_n)=(e^u)^{-\nu} f(\theta'_1,\dots,\theta'_{n-1},\theta'_n),
\label{indN}
\end{equation}
o\`u $\theta'_i$, $1\leq i\leq n$ et $e^u$ sont donn\'es par les \'egalit\'es (\ref{thet'exp1}) et (\ref{thet'exp2}).

\subsection{S\'erie sph\'erique nonunitaire pour l'alg\`ebre de Lie $\so(n+1,1)$}
Dans cette partie, nous d\'ecrivons la repr\'esentation infinit\'esimale associ\'ee \`a la s\'erie sph\'erique nonunitaire. 
Notons $dI_\nu$ (ou encore $d{\rm Ind}^G_{MAN}(0\otimes\nu)$) cette repr\'esentation, qui est donc bien entendu une repr\'esentation de l'alg\`ebre de 
Lie $\so(n+1,1)$.

\medskip

\noindent Consid\'erons tout d'abord l'\'el\'ement $H=E_{n+1,n+2}+E_{n+2,n+1}$. Nous avons
$$
h_{n+1}(t)=\exp(tH),
$$
o\`u $h_{n+1}(t)$ est donn\'e par (\ref{A}), si bien que
$$
dI_\nu(H)f(\theta_1,\dots,\theta_n)=\displaystyle\frac{d}{dt}\Big(I_\nu(h_{n+1}(t))f(\theta_1,\dots,\theta_n)\Big)|_{t=0}.
$$
A l'aide de (\ref{indA}) nous obtenons directement 
\begin{equation}
dI_\nu(H)f=\nu\cos\theta_nf~+
~\sin\theta_n\displaystyle\frac{\partial f}{\partial\theta_n}.
\label{infindh}
\end{equation}

\medskip

Nous consid\'erons \`a pr\'esent les \'el\'ements de l'alg\`ebre de Lie nilpotente $\mathfrak{n}\subset\so(n+1,1)$, 
$$
{\sf n}_i=E_{i,n+1}-E_{n+1,i}-E_{i,n+2}-E_{n+2,i}.
$$
\noindent Ces \'el\'ements constituent une base de $\mathfrak{n}$, et nous avons, pour $n({\bf a})$ d\'efini par (\ref{N}),
$$
n({\bf a})=\exp(a_1{\sf n}_1+\cdots a_n{\sf n}_n),
$$
de telle sorte que
$$
\begin{array}{rcl}
dI_\nu({\sf n}_1)f(\theta_1,\dots,\theta_n)&=&\displaystyle\frac{d}{dt}\Big(I_\nu(n(a_1,0,\dots,0))f(\theta_1,\dots,\theta_n)\Big)|_{a_1=0}\\[10pt]
dI_\nu({\sf n}_2)f(\theta_1,\dots,\theta_n)&=&\displaystyle\frac{d}{dt}\Big(I_\nu(n(0,a_2,0,\dots,0))f(\theta_1,\dots,\theta_n)\Big)|_{a_2=0}\\[10pt]
&\vdots &\\
dI_\nu({\sf n}_n)f(\theta_1,\dots,\theta_n)&=&\displaystyle\frac{d}{dt}\Big(I_\nu(n(0,0,\dots,0,a_n))f(\theta_1,\dots,\theta_n)\Big)|_{a_n=0}.\\[10pt]
\end{array}
$$
\noindent Au moyen des expressions (\ref{thet'exp1}), (\ref{thet'exp2}) et (\ref{indN}), nous aboutissons aux r\'esultats suivants:
\begin{equation}
\begin{array}{rcl}
dI_\nu({\sf n}_1)f&=&-\nu\sin\theta_1\cdots\sin\theta_n f\\[10pt]
&&+(1-\cos\theta_n)\sum_{i=1}^{n-1}~\displaystyle\frac{\sin\theta_1\cdots\sin\theta_i}{\sin\theta_i\cdots\sin\theta_n}\cos\theta_i
\displaystyle\frac{\partial f}{\partial\theta_i}\\[10pt]
&&-(1-\cos\theta_n)\sin\theta_1\cdots\sin\theta_{n-1}\displaystyle\frac{\partial f}{\partial\theta_n};\\
\end{array}
 \label{infinduced2}
\end{equation}
et, pour $1\leq i \leq n-1$,
\begin{equation}
\begin{array}{rcl}
dI_\nu({\sf n}_{i+1})f&=&
-\nu\cos\theta_i\sin\theta_{i+1}\cdots\sin\theta_n f\\[10pt]
&&-(1-\cos\theta_n)\displaystyle\frac{\sin\theta_i}{\sin\theta_{i+1}\cdots\sin\theta_n}\displaystyle\frac{\partial f}{\partial\theta_i}\\[10pt]
&&
+(1-\cos\theta_n)\sum_{j=i+1}^{n-1}~\displaystyle\frac{\sin\theta_{i+1}\cdots\sin\theta_j}{\sin\theta_j\cdots\sin\theta_n}\cos\theta_i\cos\theta_j
\displaystyle\frac{\partial f}{\partial\theta_j}\\[10pt]
&&
-(1-\cos\theta_n)\cos\theta_i\sin\theta_{i+1}\cdots\sin\theta_{n-1}\displaystyle\frac{\partial f}{\partial\theta_n}.
\end{array}
 \label{infinduced3}
\end{equation}
\begin{rema}
Dans ces formules, les notations du type $\displaystyle\frac{\sin\theta_1\cdots\sin\theta_i}{\sin\theta_i\cdots\sin\theta_n}$ ou 
$\displaystyle\frac{\sin\theta_{i+1}\cdots\sin\theta_j}{\sin\theta_j\cdots\sin\theta_n}$, bien que peu \'el\'egantes, 
permettent une meilleure lisibilit\'e pour certaines valeurs de $i$ ou de $j$. Il faut d'autre part noter que pour le cas
$i=n-1$, l'expression $\cos\theta_i\sin\theta_{i+1}\dots\sin\theta_{n-1}$ se lit $\cos\theta_{n-1}$. 
Nous adopterons cette convention dans toute la suite.
\end{rema}

\medskip

Nous avons donc enti\`erement d\'ecrit les repr\'esentations de la s\'erie sph\'erique nonunitaire de l'alg\`ebre de Lie $\so(n+1,1)$. Notre m\'ethode 
consiste 
\`a identifier $\cF_\l(\mathbb{S}^n)$ comme un module particulier de cette s\'erie. Nous allons \`a pr\'esent proc\'eder \`a cette identification.
\section{Preuve du th\'eor\`eme \ref{equiv}}
\subsection{D\'eriv\'ee de Lie de degr\'e $\l$ en coordonn\'ees sph\'eriques}
Effectuons les changements de carte n\'ecessaires pour \'etablir l'\'equivalence cherch\'ee.
\subsubsection{Changement de carte dans $\so(n+1,1)$}
En vue de notre d\'emonstration, il nous faut \'ecrire les champs de vecteurs (\ref{confGenerators}) de l'alg\`ebre de Lie conforme en termes de 
coordonn\'ees sph\'eriques, 
et en d\'eduire les expressions 
$$L_X^\l \Big(f(\theta_1,\dots,\theta_n)(d\theta_1\wedge\cdots\wedge d\theta_n)^\l\Big).$$

\noindent Nous n'aurons en fait besoin que des expressions de $X_0$ et des champs non affines $\bar X_i$, pour $1\leq i\leq n$; en effet des 
v\'erifications directes montrent que 
 le champ $X_0$ engendre la sous-alg\`ebre $\mathfrak{a}$ de dimension $1$, et nous avons, avec la repr\'esentation adjointe dans $\so(n+1,1)$,
$$
ad_{X_0}(X_i)=-X_i,\quad ad_{X_0}(\bar X_i)=+\bar X_i,
$$
si bien que l'alg\`ebre de Lie nilpotente $\mathfrak{n}=\sum_{\gamma>0}\mathfrak{g}_\gamma$ a pour g\'en\'erateurs les champs 
$\bar X_i$, $1\leq i\leq n$.

\medskip

Les coordonn\'ees st\'er\'eographiques et sph\'eriques sont li\'ees par les relations
$$
s_n=\displaystyle\frac{\sin\theta_1\cdots\sin\theta_n}{1-\cos\theta_n},\quad s_{n-i}=
\displaystyle\frac{\cos\theta_i\sin\theta_{i+1}\cdots\sin\theta_n}{1-\cos\theta_n},1\le i\le n-1.
$$
Nous obtenons alors par des calculs directs
\begin{itemize}
\item
$X_0~=~-\sin\theta_n\displaystyle\frac{\partial}{\partial\theta_n}$;
\item pour $1\leq i\leq {n-1}$,
$$
\begin{array}{rcl}
\bar X_{n-i}&=&-(1+\cos\theta_n)\displaystyle\frac{\sin\theta_i}{\sin\theta_{i+1}\sin\theta_n}\displaystyle\frac{\partial}{\partial\theta_i}\\[10pt]
&&+(1+\cos\theta_n)\sum_{j=i+1}^{n-1}\displaystyle\frac{\sin\theta_1\cdots\sin\theta_j}{\sin\theta_1\cdots\sin\theta_i\sin\theta_j\cdots\sin\theta_n}
\cos\theta_i\cos\theta_j\displaystyle\frac{\partial}{\partial\theta_j}\\[10pt]
&&+(1+\cos\theta_n)\displaystyle\frac{\sin\theta_i\cdots\sin\theta_{n-1}}{\sin\theta_i}\cos\theta_i\displaystyle\frac{\partial}{\partial\theta_n};
\end{array}
$$
\item et enfin,
$$
\begin{array}{rcl}
\bar X_{n}&=&(1+\cos\theta_n)\sum_{i=1}^{n-1}\displaystyle\frac{\sin\theta_1\cdots\sin\theta_i}{\sin\theta_i\cdots\sin\theta_n}\cos\theta_i
\displaystyle\frac{\partial}{\partial\theta_i}\\[10pt]
&&+(1+\cos\theta_n)\sin\theta_1\cdots\sin\theta_{n-1}\displaystyle\frac{\partial}{\partial\theta_n}.
\end{array}
$$
\end{itemize}
\subsubsection{D\'eriv\'ee de Lie de degr\'e $\l$}

Nous pouvons \`a pr\'esent d\'ecrire les actions $L_{X_0}^\l$ et $L_{\bar X_i}^\l$, $1\leq i\leq n$ en termes de coordonn\'ees sph\'eriques: nous avons, 
en appliquant (\ref{LieDer}),
$$
\begin{array}{rcl}
L_{X_0}^\l &=&-\sin\theta_n\displaystyle\frac{\partial }{\partial\theta_n}-\l\cos\theta_n ,
\end{array}
$$
l'op\'erateur $\cos\theta_n$ du membre  de droite \'etant un simple op\'erateur de multiplication.
De m\^eme, pour $1\leq i\leq n-1$, 
$$
\begin{array}{rcl}
L_{\bar X_{n-i}}^\l &=&-(1+\cos\theta_n)\displaystyle\frac{\sin\theta_i}{\sin\theta_{i+1}\sin\theta_n}\displaystyle\frac{\partial }{\partial\theta_i}\\[10pt]
&&+(1+\cos\theta_n)\sum_{j=i+1}^{n-1}\displaystyle\frac{\sin\theta_1\cdots\sin\theta_j}{\sin\theta_1\cdots\sin\theta_i\sin\theta_j\cdots\sin\theta_n}
\cos\theta_i\cos\theta_j\displaystyle\frac{\partial }{\partial\theta_j}\\[10pt]
&&+(1+\cos\theta_n)\displaystyle\frac{\sin\theta_i\cdots\sin\theta_{n-1}}{\sin\theta_i}\cos\theta_i\displaystyle\frac{\partial }{\partial\theta_n}\\[10pt]
&&+\l\Big(-(1+\cos\theta_n)\displaystyle\frac{\cos\theta_i}{\sin\theta_{i+1}\cdots\sin\theta_n}\\[10pt]
&&\hspace{0.8cm} -(1+\cos\theta_n)\sum_{j=i+1}^{n-1}\displaystyle\frac{\sin\theta_1\cdots\sin\theta_j}
{\sin\theta_1\cdots\sin\theta_i\sin\theta_j\cdots\sin\theta_n}
\cos\theta_i\sin\theta_j \\[10pt]
&&\hspace{0.8cm}-\sin\theta_n\cos\theta_i\displaystyle\frac{\sin\theta_i\cdots\sin\theta_{n-1}}{\sin\theta_i}
 \Big);
\end{array}
$$
et enfin :
$$
\begin{array}{rcl}
L_{\bar X_{n}}^\l &=&(1+\cos\theta_n)\sum_{i=1}^{n-1}\displaystyle\frac{\sin\theta_1\cdots\sin\theta_i}{\sin\theta_i\cdots\sin\theta_n}\cos\theta_i
\displaystyle\frac{\partial }{\partial\theta_i}\\[10pt]
&&+(1+\cos\theta_n)\sin\theta_1\cdots\sin\theta_{n-1}\displaystyle\frac{\partial }{\partial\theta_n}\\[10pt]
&&+\l\Big(-(1+\cos\theta_n)\sum_{i=1}^{n-1}\displaystyle\frac{\sin\theta_1\cdots\sin\theta_i}{\sin\theta_i\cdots\sin\theta_n}\sin\theta_i
 \\[10pt]
&&\hspace{0.8cm}-\sin\theta_1\cdots\sin\theta_n\Big).
\end{array}
$$
Nous sommes \`a pr\'esent \`a m\^eme de d\'emontrer l'identification de la d\'eriv\'ee de Lie de degr\'e $\l$ \`a 
la repr\'esentation infinit\'esimale associ\'ee \`a une repr\'esentation particuli\`ere de la s\'erie sph\'erique nonunitaire du groupe $\SO_0(n+1,1)$, 
que nous allons d\'eterminer.
\begin{defi}
Nous d\'esignons par ${\mathcal C}^\infty_\nu(\mathbb{S}^n)$ le sous-module de ${\mathcal L}^2(\mathbb{S}^n$) form\'e des \'el\'ements de classe 
${\mathcal C}^\infty$, muni de la repr\'esentation 
$d{\rm Ind}_{MAN}^G(0\otimes \nu)=dI_\nu$ de $\mathfrak{g}$, 
et de la repr\'esentation r\'eguli\`ere \`a gauche de $K$, donn\'ee par (\ref{reg}).
\end{defi}
\noindent Nous nous proposons ici de d\'emontrer l'\'equivalence des $\mathfrak{g}$- et $K$-modules $\cF_\l(\mathbb{S}^n)$ et 
${\mathcal C}^\infty_{n\l}(\mathbb{S}^n)$.
\subsection{Un r\'esultat de cohomologie}
Pour les d\'efinitions cohomologiques, on pourra consulter l'ouvrage \cite{FU} de M.~Fuks.

Nous donnons avant tout un r\'esultat g\'en\'eral, qui affirme que deux repr\'esentations de $\mathfrak{g}$ qui diff\`erent de la d\'eriv\'ee de Lie par 
deux cocycles cohomologues sont \'equivalentes
 ($L_X$ d\'esigne la d\'eriv\'ee de Lie le long du champ $X$).
\begin{lemm}
Soient $c_1$ et $c_2\in {\sf H}^1(\so(n+1,1),{\mathcal C}^\infty(\bbR^n))$. 
Nous avons
$$
\bar c_1-\bar c_2=\bar 0\implies L_X+c_1\cong L_X+c_2 .
$$
\label{cohom}
\end{lemm}
En effet, par hypoth\`ese il existe $\varphi\in{\mathcal C}^\infty(\bbR^n)$ tel que $c_1-c_2={\rm d}\varphi$. Alors l'application 
$$
\begin{array}{rccl}
\psi :&\cF_\l(\mathbb{S}^n)&\longrightarrow&\cF_\l(\mathbb{S}^n)\\
&f&\longmapsto&{\rm e}^\varphi f
\end{array}
$$
entrelace les actions $L_X+c_1$ et $L_X+c_2$: nous avons $c_1(X)-c_2(X)=L_X\varphi$ pour tout $X$, de sorte que
$$
\begin{array}{rcl}
L_X+c_2(X)({\rm e}^{\varphi}f)&=&{\rm e}^{\varphi}L_X f+L_X(\varphi) f.{\rm e}^{\varphi}+c_2(X)f{\rm e}^{\varphi}\\
&=&{\rm e}^{\varphi}\Big(L_X f+L_X(\varphi)f+c_2(X)\Big)\\
&=&{\rm e}^{\varphi}\Big(L_X f+c_1(X)\Big),
\end{array}
$$
ce qui d\'emontre le lemme \ref{cohom}.

\medskip

Notre d\'emonstration de l'\'equivalence est bas\'ee sur deux \'equivalences successives.
\subsection{Premi\`ere \'etape}
Notons tout d'abord $M^\nu$ la repr\'esentation de l'alg\`ebre de Lie conforme sur $\cF_\l(\mathbb{S}^n)$ donn\'ee en coordonn\'ees sph\'eriques par 
l'action
$$
M^\nu _X=L_X+\nu\displaystyle\frac{\partial X^n}{\partial\theta_n}
$$
(on notera que $\l$ n'intervient pas dans cette action).
Une premi\`ere \'etape dans notre d\'emarche est donn\'ee par le lemme suivant :
\begin{lemm}
Les repr\'esentations $dI_\nu$ et $M^\nu$ sont \'equivalentes.
\end{lemm}
\noindent En effet nous identifions les bases vectorielle et matricielle de $\so(n+1,1)$ par la correspondance $X_0\mapsto H$, 
$\bar X_n\mapsto -{\sf n}_1$ 
et $\bar X_{n-i}\mapsto -{\sf n}_{i+1},~1\le i\le n-1$, 
si bien qu'\`a l'aide des formules (\ref{infindh}), (\ref{infinduced2}) et (\ref{infinduced3}), nous obtenons directement pour tout $X\in\so(n+1,1)$
$$
\pi\circ dI_\nu(X)=M^\nu_X\circ\pi,
$$
o\`u
$$
\begin{array}{rccl}
\pi:&{\mathcal C}^\infty_\nu(\mathbb{S}^n)&\longrightarrow&\cF_\l(\mathbb{S}^n)\\
&f(\theta_1,\dots,\theta_n)&\longmapsto&f(\theta_1,\dots,\theta_{n-1},\pi+\theta_n)(d\theta_1\wedge\cdots\wedge d\theta_n)^\l.
\end{array}
$$

Il nous reste \`a pr\'esent \`a d\'emontrer que les repr\'esentations $M^\nu_X$ et $L_X^\l$, qui s'effectuent toutes deux sur $\cF_\l(\mathbb{S}^n)$,
 sont \'equivalentes si et seulement si $\nu=n\l$. Nous nous appuierons \`a cet effet sur le lemme \ref{cohom}.
\subsection{Deuxi\`eme \'etape et \'equivalence}
Avant tout, constatons les deux faits suivants, v\'erifi\'es pour tout $X=\sum_i X^i\displaystyle\frac{\partial}{\partial\theta_i}$ de l'alg\`ebre de Lie 
$\so(n+1,1)$ : 

- $L_X^\l$ est du type $L_X+c_1(X)$, avec $c_1(X)=\l\Div X$, o\`u $\Div X=\sum_i\displaystyle\frac{\partial X^i}{\partial\theta_i}$;

- $M_X^\nu$ est du type $L_X+c_2(X)$, avec $c_2(X)=\nu\displaystyle\frac{\partial X^n}{\partial\theta_n}$.

\noindent Nous avons donc avec ces notations 
$$
c_1(X)-c_2(X)=\l(\Div X-\displaystyle\frac{\nu}{\lambda}\displaystyle\frac{\partial X^n}{\partial\theta_n}).
$$
 Nous allons \`a pr\'esent prouver que $\Div X-\displaystyle\frac{\nu}{\lambda}\displaystyle\frac{\partial X^n}{\partial\theta_n}$
 est un cobord dans l'espace de cohomologie ${\sf H}^1(\so(n+1,1),{\mathcal C}^\infty(\bbR^n))$. Nous prouvons pr\'ecis\'ement que c'est le cas si et 
seulement si $\nu=n\l$.
\begin{prop}
$$
\Div X-\alpha\displaystyle\frac{\partial X^n}{\partial\theta_n}\in\bar 0\in{\sf H}^1(\so(n+1,1),{\mathcal C}^\infty(\bbR^n))\iff \alpha=n.
$$
\label{cob}
\end{prop}
\noindent En effet, si nous \'ecrivons l'\'equation 
\begin{equation}
\Div(X) -\alpha\displaystyle\frac{\partial X^n}{\partial\theta_n}=d\varphi(X)
\label{sys}
\end{equation}
 pour chaque \'el\'ement $X$ de $\so(n+1,1)$, nous
 obtenons les $n+1$ \'equations suivantes :

\begin{itemize}
\item pour $X=\bar X_n$,
$$
\begin{array}{l}
(\alpha-1)\sin\theta_1\cdots\sin\theta_n-(1+\cos\theta_n)\sum_{i=1}^{n-1}\displaystyle\frac{\sin\theta_1\cdots\sin\theta_i}{\sin\theta_i\cdots\sin\theta_n}
\sin\theta_i\\[10pt]
=(1+\cos\theta_n)\Big(\sum_{i=1}^{n-1}\displaystyle\frac{\sin\theta_1\cdots\sin\theta_i}{\sin\theta_i\cdots\sin\theta_n}
\cos\theta_i\displaystyle\frac{\partial \varphi}{\partial\theta_i}
+\sin\theta_1\cdots\sin\theta_{n-1}\displaystyle\frac{\partial \varphi}{\partial\theta_n}\Big);
\end{array}
$$
\item pour les $\bar X_{n-i}$, $1\le i\le n-1$, nous avons $n-1$ \'equations index\'ees par $i$ :
\begin{equation}
\begin{array}{l}
(\alpha-1)\cos\theta_i\sin\theta_{i+1}\cdots\sin\theta_n\\[10pt]-(1+\cos\theta_n)\Big(\displaystyle\frac{\cos\theta_i}{\sin\theta_{i+1}\cdots\sin\theta_n}\\[10pt]
\hspace{60pt}+
\sum_{j=i+1}^{n-1}\displaystyle\frac{\sin\theta_1\cdots\sin\theta_j}{\sin\theta_1\cdots\sin\theta_i\sin\theta_j\cdots\sin\theta_n}
\cos\theta_i\sin\theta_j\Big)\\[10pt]
=(1+\cos\theta_n)\Big(-\displaystyle\frac{\sin\theta_i}{\sin\theta_{i+1}\cdots\sin\theta_n}\displaystyle\frac{\partial \varphi}{\partial\theta_i}\\[10pt]
\hspace{70pt}+\sum_{j=i+1}^{n-1}\displaystyle\frac{\sin\theta_1\cdots\sin\theta_j}{\sin\theta_1\cdots\sin\theta_i\sin\theta_j\cdots\sin\theta_n}
\cos\theta_i\cos\theta_j\displaystyle\frac{\partial \varphi}{\partial\theta_j}\\[10pt]
\hspace{70pt}+\cos\theta_i\sin\theta_{i+1}\cdots\sin\theta_{n-1}\displaystyle\frac{\partial \varphi}{\partial\theta_n}\Big);
\end{array}
\label{eqi}
\end{equation}
\item enfin, pour $X=X_0$, 
$$
(\alpha -1)\cos\theta_n=-\sin\theta_n\displaystyle\frac{\partial \varphi}{\partial\theta_n}.
$$
\end{itemize}
Le syst\`eme lin\'eaire form\'e de ces $n+1$ \'equations admet une unique solution donn\'ee par 
$$
\left\{
\begin{array}{rcl}
\displaystyle\frac{\partial \varphi}{\partial\theta_k}&=&-\displaystyle\frac{k-1}{\tan\theta_k},~k=1,\dots,n\\[10pt]
\alpha&=&n.
\end{array}
\right.
$$

\medskip 

\noindent En effet, pour r\'esoudre ce syst\`eme nous pouvons proc\'eder de la fa\c{c}on suivante :
\begin{itemize}
\item Appelons $S_1$ le syst\`eme donn\'e ci-dessus, et num\'erotons de $(1)$ \`a $(n+1)$ les \'equations, dans l'ordre d\'ecrit: ainsi, l'\'equation $(1)$ 
est obtenue 
en appliquant l'\'equation ($\ref{sys}$) au champ $\bar X_n$, et pour $i$ variant de $2$ \`a $n$, l'\'equation $(i)$ est obtenue en 
appliquant cette \'equation ($\ref{sys}$) au champ $\bar X_{n-i+1}$. Enfin, l'\'equation $(n+1)$ est obtenue par ce m\^eme proc\'ed\'e, au moyen du 
champ $ X_0$.
\begin{enumerate}
\item Multiplions l'\'equation $(1)$ par $\cos\theta_1$, et soustrayons au r\'esultat obtenu le produit de l'\'equation $(2)$ par $\sin\theta_1$: nous 
obtenons $\displaystyle\frac{\partial \varphi}{\partial\theta_1}=0$. 
\item Multiplions l'\'equation $(1)$ par $\sin\theta_1$, et additionnons au r\'esultat obtenu le produit de l'\'equation $(2)$ par $\cos\theta_1$: nous 
obtenons une nouvelle \'equation, o\`u ne figurent plus les termes o\`u apparaissait l'angle $\theta_1$. 
Nous construisons alors un syst\`eme $S_2$, compos\'e de $n$ \'equations num\'erot\'ees de $(2)$ \`a $(n+1)$, o\`u cette nouvelle \'equation porte le 
num\'ero $(2)$, et les autres \'equations sont inchang\'ees par rapport au syst\`eme $S_1$.
\end{enumerate}
\item R\'ep\'etons le proc\'ed\'e d\'ecrit ci-dessus, permettant, pour $i$ de $(1)$ \`a $(n-1)$, de passer du syst\`eme $S_i$ au syst\`eme $S_{i+1}$: 
dans $S_i$, les \'equations sont num\'erot\'ees de $(i)$ \`a $(n+1)$. On proc\`ede en deux temps :
\begin{enumerate}
\item[a.] $(i)\cos\theta_i-(i+1)\sin\theta_i\iff\displaystyle\frac{\partial \varphi}{\partial\theta_i}=-\displaystyle\frac{i-1}{\tan\theta_i}$ ;
\item[b.] $(i)\sin\theta_i+(i+1)\cos\theta_i$ est l'\'equation $(i+1)$ du nouveau syst\`eme $S_{i+1}$.
\label{algo}
\end{enumerate}
Explicitement, l'\'equation $(i)$ du syst\`eme $S_i$ est la suivante :
$$\begin{array}{l}
(\alpha-1)\sin\theta_i\cdots\sin\theta_n-(1+\cos\theta_n)\Big(\displaystyle\frac{\sin\theta_i}{\sin\theta_{i+1}\cdots\sin\theta_n}+\displaystyle\frac{i-1}
{\sin\theta_i\cdots\sin\theta_n}
\\[10pt]\hspace{100pt}+\sum_{j=i+1}^{n-1}\displaystyle\frac{\sin\theta_1\cdots\sin\theta_j}{\sin\theta_1\cdots\sin\theta_i\sin\theta_j\cdots\sin\theta_n}
\sin\theta_i\sin\theta_j\Big)\\[15pt]
=(1+\cos\theta_n)\Big(\displaystyle\frac{\cos\theta_i}{\sin\theta_{i+1}\cdots\sin\theta_n}\displaystyle\frac{\partial \varphi}{\partial\theta_i}\\[10pt]
\hspace{70pt}+ \sum_{j=i+1}^{n-1}\displaystyle\frac{\sin\theta_1\cdots\sin\theta_j}{\sin\theta_1\cdots\sin\theta_i\sin\theta_j
\cdots\sin\theta_n}\sin\theta_i\cos\theta_j
\displaystyle\frac{\partial \varphi}{\partial\theta_j}\\[10pt]\hspace{70pt}+\sin\theta_i\cdots\sin\theta_{n-1}
\displaystyle\frac{\partial \varphi}{\partial\theta_n}\Big),
\end{array}$$
et l'\'equation $(i+1)$ du syst\`eme $S_i$ est simplement (\ref{eqi}).

Nous v\'erifions alors facilement par r\'ecurrence les assertions a. et b. ci-dessus.
\item Le syst\`eme $S_n$ s'\'ecrit enfin

$$
\begin{array}{rcl}
(\alpha-1)\sin\theta_n-(1+\cos\theta_n)\displaystyle\frac{n-1}{\sin\theta_n}&=&(1+\cos\theta_n)\displaystyle\frac{\partial \varphi}{\partial\theta_n}\\[10pt]
(\alpha-1)\cos\theta_n&=&-\sin\theta_n\displaystyle\frac{\partial \varphi}{\partial\theta_n}
\end{array}
$$
et ce dernier syst\`eme admet clairement une solution unique, donn\'ee par
$$\displaystyle\frac{\partial \varphi}{\partial\theta_n}=-\displaystyle\frac{n-1}{\tan\theta_k} \quad {\rm et}\quad  \alpha=n.$$ 
\end{itemize}
Le syst\`eme $S_1$ est ainsi r\'esolu, et la proposition \ref{cob} est d\'emontr\'ee.

\bigskip

Ainsi en appliquant le lemme \ref{cohom}, nous avons d\'emontr\'e que les $\mathfrak{g}$-modules $\cF_\l(\mathbb{S}^n)$ et 
${\mathcal C}^\infty_{n\l}(\mathbb{S}^n)$ sont isomorphes. 
Il est de plus clair que cet isomorphisme est compatible avec l'action de $K$, puisque l'action du groupe $K$ est donn\'ee sur chaque 
espace par la m\^eme repr\'esentation r\'eguli\`ere \`a gauche. Ainsi, si nous d\'esignons par $\varphi$ l'isomorphisme de 
$\mathfrak{g}$-modules de $\cF_\l(\mathbb{S}^n)$ vers ${\mathcal C}^\infty_\nu(\mathbb{S}^n)$, il vient imm\'ediatement
$$
\varphi(k.(f(\theta_1,\dots,\theta_n)(d\theta_1\wedge\cdots\wedge d\theta_n)^\l))=k.\varphi(f(\theta_1,\dots,\theta_n)
d\theta_1\wedge\cdots\wedge d\theta_n)^\l),
$$
pour tout $k\in K$.

Le th\'eor\`eme \ref{equiv} est donc d\'emontr\'e.
Nous \'ecrivons
$$
L_X^\l\cong d{\rm Ind}_{MAN}^G(0\otimes n\l)(X).
$$

 \medskip

Nous d\'ecrivons \`a pr\'esent en d\'etail les $({\mathfrak g},K)$-modules de la s\'erie sph\'erique nonunitaire, dans le but de leur identifier ceux de 
$\cF_\l(\mathbb{S}^n)$.
\section {$({\mathfrak g},K)$-modules simples de la s\'erie sph\'erique nonunitaire}
Nous devons avant tout donner des indications pr\'ecises sur les repr\'esentations du groupe $K=\SO(n+1)$. Deux cas sont \`a consid\'erer, selon 
la parit\'e de 
$n+1$. Nous ne donnons les d\'etails que pour $n+1=2k$, le travail \'etant similaire pour $n+1$ impair.

Les notations utilis\'ees dans tout ce qui suit sont celles utilis\'ees par M.~Guichardet dans (\cite{Gui}, appendice $B.10$).
\subsection{Poids et repr\'esentations de $\SO(n+1)$}
Nous posons ici $n+1=2k$, c'est \`a dire $k=[\frac{n+1}{2}]$, o\`u $[p]$ est la partie enti\`ere de $p$.

\medskip

\noindent Rappelons avant tout la d\'efinition des racines et poids :
Soit ${\mathfrak g}$ une alg\`ebre de Lie semi-simple, et ${\mathfrak h}$ une sous-alg\`ebre de Cartan de ${\mathfrak g}$.

- On dit que $\alpha\in{\mathfrak h}^*$ est une racine de $({\mathfrak g},{\mathfrak h})$ si le sous-espace ${\mathfrak g}_\alpha$ d\'efini par
$$
{\mathfrak g}_\alpha=\{X\in{\mathfrak g}~:~[H,X]=\alpha(H)X~\forall H\in{\mathfrak h}\}
$$
est non r\'eduit \`a $\{0\}$.

- On dit que $\l\in{\mathfrak h}^*$ est un poids d'une repr\'esentation $\tau$ de ${\mathfrak g}$ sur un espace $V$ si le sous-espace $V_\l$ d\'efini par
$$
V_\l=\{v\in V~:~\tau(H)v=\l(H)v~\forall H\in{\mathfrak h}\}
$$
est non r\'eduit \`a $\{0\}$.
Ainsi, une racine de $({\mathfrak g},{\mathfrak h})$ est un poids pour la repr\'esentation adjointe dans ${\mathfrak g}$.

\bigskip

En ce qui nous concerne, prenons d'abord une base de l'alg\`ebre $\mathfrak{k}=\so(n+1)$ et notamment d'une sous alg\`ebre de Cartan :
cette sous-alg\`ebre de Cartan $\mathfrak{h}$ est engendr\'ee par les \'el\'ements 
$$
H_{r}=i(E_{2r-1,2r}-E_{2r,2r-1}),\quad r=1,2,\dots, k.
$$
Consid\'erons les applications lin\'eaires sur cette sous-alg\`ebre de Cartan donn\'ees par
$$
\varepsilon_r(\l_1H_1+\cdots+\l_kH_k)=\l_r.
$$
Les racines de $(\so(n+1),\mathfrak{h})$ sont, pour $1\leq r<s\leq k$, les applications
$$
\varepsilon_r-\varepsilon_s,~\varepsilon_r+\varepsilon_s,~-\varepsilon_r+\varepsilon_s,~-\varepsilon_r-\varepsilon_s.
$$
Nous ne donnons pas les \'el\'ements et sous-espaces radiciels, inutiles ici.

- Les racines positives sont $\varepsilon_r-\varepsilon_s,~\varepsilon_r+\varepsilon_s$.

- La plus grande racine est $\varepsilon_1+\varepsilon_2$ (poids dominant de la repr\'esentation adjointe.)

- Les poids fondamentaux ($\mathbb{Z}$-base du r\'eseau des poids) sont 
$$
\begin{array}{rcl}
\varpi_r&=&\varepsilon_1+\varepsilon_2+\cdots+\varepsilon_r,~r=1,\dots,k-2 ;\\[5pt]
\varpi_{k-1}&=&\frac{1}{2}(\varepsilon_1+\varepsilon_2+\cdots+\varepsilon_{k-1}-\varepsilon_k) ;\\[5pt]
\varpi_k&=&\frac{1}{2}(\varepsilon_1+\varepsilon_2+\cdots+\varepsilon_{k-1}+\varepsilon_k).
\end{array}
$$
Par cons\'equent si un poids $\Lambda$ s'\'ecrit 
$$
\Lambda=\sum_{i=1}^{k}\l_i\varpi_i=\sum_{i=1}^{k}\m_i\varepsilon_i, \l_i\in\mathbb{Z},
$$
nous aurons 
$$
\l_1=\m_1-\m_2,~\l_2=\m_2-\m_3,\dots,~\l_{k-1}=\m_{k-1}-\m_k,~\l_k=\m_{k-1}+\m_k.
$$
$\Lambda$ est un poids dominant si et seulement si tous les $\l_i$ sont positifs, ce qui nous donne 
$$
\m_1\geq\m_2\geq\cdots\geq\m_{k-1}\geq |\m_k|.
$$
Les $\m_i$ sont entiers ou demi-entiers non entiers. Ces poids dominants sont aussi ceux du groupe ${\rm Spin}(2k)$, rev\^etement universel \`a deux 
feuillets de $\SO(2k)$.
Ceux de $\SO(2k)$ sont obtenus avec les $\m_i$ tous entiers.
 
Une repr\'esentation est caract\'eris\'ee par son plus haut poids, et ce dernier est forc\'ement un poids dominant.
Ainsi une repr\'esentation de $\SO(2k)$ est not\'ee 
$$
D^{m_1,\dots,m_k},~ {\rm avec}~ |m_1|\leq m_2\leq\cdots\leq m_{k-1}\leq m_k,
$$
o\`u nous avons not\'e $m_1=\m_k$, $m_2=\m_{k-1},\dots m_k=\m_1$.

\medskip

De m\^eme, dans le cas o\`u $n+1$ est impair, $n+1=2k+1$, on d\'emontre que pour une repr\'esentation de $\SO(2k+1)$ on peut noter
$$
D^{m_1,\dots,m_k},\quad 0\leq m_1 \leq m_2\leq\cdots\leq m_{k-1}\leq m_k.
$$
\subsection{$(\mathfrak{g},K)$-modules} 
Nous donnons tout d'abord quelques indications concernant les $(\mathfrak{g},K)$-modules, dont on pourra trouver les d\'etails dans l'ouvrage 
\cite{Kna} de M.~Knapp par exemple.
 
\noindent Il convient tout d'abord de remarquer que $\cF_\l(\mathbb{S}^n)$ satisfait aux conditions de d\'efinition d'un $(\mathfrak{g},K)$-module 
donn\'ees en introduction de cet article, si ce n'est que tout vecteur de cet espace n'est pas forc\'ement $K$-fini : en effet, nous avons une 
repr\'esentation (analytique) que nous notons $\pi$ de $G=\SO_0(n+1,1)\subset \Diff(\mathbb{S}^n$). 
D'autre part dans $G$ nous avons $$k{\rm exp}(tX)k^{-1}={\rm exp}({\rm Ad}k\cdot tX)$$ pour tous $k\in K$, $t\in\bbR$, $X\in\mathfrak{g}$ 
(voir par exemple \cite{Bou}, par. $6$, cor. $2$). En prenant la d\'eriv\'ee en $t=0$ de 
$\pi(k{\rm exp}(tX)k^{-1})\cdot v$, la diff\'erentielle 
$d\pi:\mathfrak{g}\longrightarrow\End (\cF_\l(\mathbb{S}^n))$ est telle que l'on obtient $$d\pi({\rm Ad}(k)\cdot X)\cdot v=k\cdot(d\pi(X)\cdot (k^{-1}\cdot v)).$$

\medskip

Lorsque $K$ agit au moyen d'op\'erateurs unitaires la repr\'esentation $\pi$ de $G$ restreinte \`a $K$ se d\'ecompose en somme directe orthogonale
 d'espaces sur lesquels elle est irr\'eductible,                                         
 et la multiplicit\'e de chaque repr\'esentation est bien d\'efinie :
\begin{equation}
\pi|_K ~\cong~\sum_{\tau\in\hat{K}}n_\tau \tau
\label{Ktypes}
\end{equation}
o\`u $\hat{K}$ est l'ensemble des classes d'\'equivalence  des repr\'esentations irr\'eductibles de $K$ et o\`u chaque multiplicit\'e $n_\tau$ est un entier 
positif ou est $+\infty$. 
Un vecteur d'un de ces irr\'eductibles $K$-espaces est $K$-fini, de sorte que l'espace des vecteurs $K$-finis est dense dans $E$. 
Les classes d'\'equivalence $\tau$ dans ($\ref{Ktypes})$ de multiplicit\'e non nulle sont appel\'ees les $K$-types de $\pi$, et $n_\tau \tau$ 
est la composante isotypique de type $\tau$.

\noindent Une telle repr\'esentation $\pi$ est dite admissible si $\pi(K)$ op\`ere par des op\'erateurs unitaires et si chaque $\tau$ dans $\hat{K}$ est 
de multiplicit\'e finie.

Nous allons voir que c'est le cas pour ce qui nous concerne: M.~Vilenkin a d\'emontr\'e dans \cite{Vil} que la repr\'esentation r\'eguli\`ere \`a gauche de 
$\SO(n+1)$ est unitaire pour le produit scalaire usuel dans ${\mathcal L}^2(\mathbb{S}^n)$
et se d\'ecompose en la somme directe 
$${\mathcal H}(K)=\bigoplus {\mathcal H}^{n+1,m},~m\in\mathbb{N}$$
des sous $K$-modules ${\mathcal H}^{n+1,m}$ des polyn\^omes harmoniques homog\`enes en les coordonn\'ees cart\'esiennes, chacun de ces 
sous-modules \'etant simple 
(except\'e pour le cas $n=1$ o\`u ${\mathcal H}^{2,m}$ est somme directe de deux sous-modules simples de dimension $1$, du fait du caract\`ere 
ab\'elien de $\SO(2)$), 
et de multiplicit\'e $1$. L'isomorphisme donn\'e par le th\'eor\`eme \ref{equiv} nous montre donc que notre repr\'esentation est admissible, \`a 
condition bien s\^ur d'avoir prouv\'e que les ${\mathcal H}^{n+1,m}$ d\'ecrivent de mani\`ere exhaustive l'espace des vecteurs $K$-finis.

\subsection{$(\mathfrak{g},K)$-modules pour ${\rm Ind}_{MAN}^G(0\otimes \nu+\rho)$}
Afin que nos notations soient en accord avec celles donn\'ees dans (\cite{Gui}, appendice B.10), nous \'etudions la repr\'esentation induite par le 
caract\`ere $a^{\rho + \nu}$ de $A$ et non 
$a^{\nu}$ comme nous l'avions fait auparavant. Ainsi, nous \'etudions les vecteurs $K$-finis de la repr\'esentation 
${\rm Ind}_{MAN}^G(0\otimes \nu+\rho)$, et non 
${\rm Ind}_{MAN}^G(0\otimes \nu)$.
Il convient de noter d\'es \`a pr\'esent que nous avons d\'emontr\'e l'\'equivalence 
$$
L_X^\l\cong d{\rm Ind}_{MAN}^G(0\otimes n\l)(X).
$$
Aussi nous faudra-t-il remplacer $\nu+\rho$ par $n\l$ dans les r\'esultats qui suivent.
 
 \bigskip

Notons $E_{0, \nu}$  le $(\mathfrak{g},K)$-module des vecteurs $K$-finis de ${\rm Ind}_{MAN}^G(0\otimes \nu+\rho)$.
Nous avons 
\begin{equation}
E_{0,\nu}|_K=\bigoplus D^{0,\dots,0,m},\left|\begin{array}{l}m\in\mathbb{N}~{\rm pour}~n>1\\m\in\mathbb{Z}~{\rm pour}~n=1\end{array}\right.
\label{Kf}
\end{equation}
\begin{rema}
Nous utilisons la propri\'et\'e fondamentale  $E^*_{0,\nu}\cong E_{0,-\nu}$ (dual $K$-fini).
\end{rema}
\subsubsection{R\'esultats pour $G=\SO_0(2,1)$}
Nous r\'ef\'erant \`a (\cite{Gui}, appendice B.10), nous posons $n_1=+\infty$, et interpr\'etons les r\'esultats donn\'es pour $n=2k$.

\noindent Suite de composition de $E_  {0,\nu}$: 

$E_{0,\nu}$ est simple ssi $\nu\not\in\frac{1}{2}+\mathbb{Z}$.
Dans le cas contraire, il admet trois sous-quotients simples not\'es $E^0_{0,\nu}$, $E^+_{0,\nu}$, $E^-_{0,\nu}$. Leurs restrictions \`a $K$ sont 
les suivantes: 

 - L'entier $m$ de $E^0_{0,\nu}$ est caract\'eris\'e par $|m|\le |\nu|-\frac{1}{2}$.

- Celui de $E^{\pm}_{0,\nu}$ par $|\nu|+\frac{1}{2}\le \pm m$.

 - $E^0_{0,\nu}$ est un sous-module de $E_{0,\nu}$ si $\nu<0$, et un quotient si $\nu>0$.

 - Les trois sous-quotients donn\'es ci-dessus sont des $(\mathfrak{g},K)$-modules simples. 

 - Les sous-modules de dimension finie de $E_{0,\nu}$ sont les $E^0_{0,\nu}$, pour $\nu<0$. 

 \noindent Les modules unitaires sont :

- $E_{0,\nu}$ pour $\nu$ imaginaire pur (s\'erie principale unitaire).

- $E_{0,\nu}$ pour $\nu\in ]0,\frac{1}{2}[$ (s\'erie compl\'ementaire).

- $E^0_{0,\frac{1}{2}}$.

- $E^\pm_{0,\nu}$ pour $\nu\in\frac{1}{2}+\mathbb{Z}$ (s\'erie discr\`ete).

\subsubsection{R\'esultats pour $G=\SO_0(n+1,1)$, $n>1$}

Utilisant les m\^emes notations, nous constatons que les cas $n+1=2k$ et $n+1=2k+1$ donnent pour ce qui nous concerne les 
m\^emes r\'esultats, qui sont les suivants :

\noindent Suite de composition de $E_  {0,\nu}$: 

- Si $n$ est pair, $E_  {0,\nu}$ est simple ssi $\nu\not\in\mathbb{Z}$ ou $\nu=\pm(0,1,\dots,\frac{n}{2}-1)$;

- Si $n$ est impair, $E_  {0,\nu}$ est simple ssi $\nu\not\in\frac{1}{2}+\mathbb{Z}$ ou $\nu=\pm(\frac{1}{2},\frac{3}{2},\dots,\frac{n}{2}-1)$;

 \noindent Dans le cas contraire, c'est \`a dire si $\nu=\pm(\frac{n}{2},\frac{n}{2}+1,\dots)$, il admet deux sous-quotients simples not\'es 
$E^0_{0,\nu}$ et $E^\times_{0,\nu}$. L'entier $m$ de $E^0_{0,\nu}$ est caract\'eris\'e par $m\le |\nu|-k+\frac{1}{2}$, et celui de $E^\times_{0,\nu}$ 
par $m\ge |\nu|-k+\frac{3}{2}$. 

$E^0_{0,\nu}$ est un sous-module de dimension finie de $E_{0,\nu}$ si $\nu<0$, et un quotient si $\nu>0$, et on obtient ainsi les sous-$(\mathfrak{g},K)$
-modules de $E_{0,\nu}$.

\noindent Les modules unitaires sont les suivants :

- $E_{0,\nu}$ pour $\nu$ imaginaire pur (s\'erie principale unitaire).

- $E_{0,\nu}$ pour $\nu\in ]0,\frac{n}{2}[$ (s\'erie compl\'ementaire).

- $E^0_{0,\frac{n}{2}}$.

De cette classification exhaustive des $(\mathfrak{g},K)$-modules simples et unitaires de ${\rm Ind}_{MAN}^G(0\otimes \nu+\rho)$,  nous allons d\'eduire 
la classification
 inh\'erente au module $\cF_\l(\mathbb{S}^n)$.
\section {Preuve du th\'eor\`eme \ref{class} }
A l'aide du th\'eor\`eme \ref{equiv}, nous pouvons identifier les vecteurs $K$-finis de $\cF_\l(\mathbb{S}^n)$ \`a ceux de 
${\mathcal C}^\infty_\nu(\mathbb{S}^n)$
, de la mani\`ere suivante :
\subsection{Vecteurs $K$-finis}
Nous utilisons les notations donn\'ees pr\'ec\'edemment. 
Posons
$$
{\mathcal H}^{n+1,m}=\big\{\displaystyle\frac{P_m(x_0,\dots,x_n)}{(x_0^2+\cdots+x_n^2)^{\frac{m}{2}+
\l\frac{n+1}{2}}}\big\}\in{\it C}^\infty_{-\l(n+1)}(\bbR^{n+1}),
$$
o\`u $P_m$ d\'ecrit l'ensemble des polyn\^omes harmoniques et homog\`enes de degr\'e $m$ en les coordonn\'ees $x_0,\dots,x_n$.

${\mathcal H}^{n+1,m}$ est consid\'er\'e comme sous-espace de $\cF_\l(\mathbb{S}^n)$ \`a l'aide du ``rel\`evement'' (\ref{rel}). 
Nous consid\'erons cet espace comme $K$-module, i.e. comme $\SO(n+1)$-module, avec la repr\'esentation r\'eguli\`ere 
$
(k_0\cdot f)(k)=f(k_0^{-1}k).
$
Le $K$-module ${\mathcal H}^{n+1,m}$ est simple pour $n>1$, et pour $n=1$, il est la somme directe de deux $K$-modules simples, 
not\'es $H_m$ et $H_{-m}$, 
engendr\'es respectivement 
par l'\'el\'ement $(x_0+{\bf i}x_1)^m(x_0^2+x_1^2)^{-\frac{m}{2}-\l}$, et par son conjugu\'e. 

\medskip

Des v\'erifications directes (voir par exemple \cite{Goo}, th\'eor\`eme $5.2.4$, ou encore \cite{Kna}, page 89) 
montrent les faits suivants :
\begin{itemize}
\item Si $n=1$, $H_{\pm m}\cong D^{\pm m}$;
\item Si $n>1$, ${\mathcal H}^{n+1,m}$ est un $K$-module de plus haut poids $m\varepsilon_1$, avec l'\'el\'ement primitif 
$$
\psi_{1m}~=~\displaystyle\frac{(x_0+{\bf i}x_1)^m}{(x_0^2+\cdots+x_n^2)^{\frac{m}{2}+\l\frac{n+1}{2}}}.
$$
Ainsi ${\mathcal H}^{n+1,m}\cong D^{0,\dots,0,m}$.
\end{itemize}
Posons 
$$
{\mathcal H}(K)=\bigoplus_m{\mathcal H}^{n+1,m}
$$
\noindent Nous savons que l'espace des vecteurs $K$-finis de la repr\'esentation induite
${\rm Ind}_{MAN}^G(0\otimes \nu+\rho)$ est donn\'e par (\ref{Kf}). Par l'\'equivalence 
d\'emontr\'ee, nous pouvons \'enon\c{c}er le r\'esultat suivant :

\begin{prop}
Le $(\mathfrak{g},K)$-module des vecteurs $K$-finis de l'espace des densit\'es tensorielles sur la sph\`ere est isomorphe \`a l'espace 
${\mathcal H}(K)$ des polyn\^omes harmoniques en les coordonn\'ees cart\'esiennes sur la sph\`ere.
\end{prop}
\subsection{Classification }
Nous achevons au moyen des r\'esultats pr\'ec\'edents la d\'emonstration du th\'eor\`eme \ref{class}: 
 \begin{itemize}
 \item
Cas $n=1$: $E^0_{0,\nu}=\bigoplus_{|m|\le |\nu|-\frac{1}{2}}H_m$. Cet espace est un sous-module simple de dimension finie pour 
$\nu\in\frac{1}{2}-\mathbb{N}$. Ainsi, si $\nu=-\frac{1}{2}-l$, $l\in\mathbb{N}$, nous avons $\nu+\rho= -l$ (pour $\SO_0(n+1,1)$, 
$\rho=\frac{n}{2}$), d'o\`u
$|m|\le |\nu|-\frac{1}{2}\iff|m|\le l$. Par cons\'equent nous obtenons le sous-$(\mathfrak{g},K)$-module simple de dimension finie qui 
est la somme directe de $K$-modules simples
$\bigoplus_{|m|\le l}H_m$ pour $\l=-l$, dans $\cF_\l(\mathbb{S}^1)$, ce qui correspond \`a l'ensemble des polyn\^omes harmoniques 
homog\`enes de degr\'e $\le l$ en les coordonn\'ees cart\'esiennes $x_0,x_1$.

De la m\^eme fa\c{c}on, les autres sous-$(\mathfrak{g},K)$-modules simples sont obtenus pour $\l=l$, $l\in\mathbb{N}^*$, et sont 
les sous-espaces donn\'es par $E^{+}_{0,\nu}|_K=\bigoplus_{l\le  m} H_m$ et $E^{-}_{0,\nu}|_K=\bigoplus_{-m\le -l } H_{-m}$. Remarquons 
que la somme directe de $E^{+}_{0,\nu}$ et $E^{-}_{0,\nu}$ est l'ensemble des polyn\^omes harmoniques de degr\'e $\ge l$.
\item Cas $n>1$: 

- Si $\nu=-\frac{n}{2}-l$, $l\in\mathbb{N}$, ce qui \'equivaut
 \`a $\nu+\rho=-l$, il existe un unique sous-$(\mathfrak{g},K)$-module simple de ${\mathcal H}(K)$, de dimension finie, donn\'e par le sous-espace  
$E^0_{0,\nu}|_K=\bigoplus_{m\le |\nu|-\frac{n}{2}}{\mathcal H}^{n+1,m}$.  Ce sous-$(\mathfrak{g},K)$-module de $\cF_\l(\mathbb{S}^n)$ est 
donc donn\'e, pour $l+n\l=0$, par 
l'ensemble des polyn\^omes harmoniques homog\`enes de degr\'e $\le l$ en les coordonn\'ees cart\'esiennes (en effet $m\le
 |\nu|-\frac{n}{2}\iff m\le -n\l$).

- Si $\nu=\frac{n}{2}+l$, $l\in\mathbb{N}$, ce qui \'equivaut \`a $\nu+\rho=l+n$, on obtient un unique sous-$(\mathfrak{g},K)$-module simple 
de dimension infinie,  qui est l'espace des 
polyn\^omes harmoniques homog\`enes de degr\'e $\ge l+1$; ce cas correspond \`a la valeur $n\l=l+n\iff \l=1+l/n$.
\end{itemize}

\medskip

Il reste \`a traduire les conditions d'unitarisabilit\'e de chacun ce ces $(\mathfrak{g},K)$-modules. Nous obtenons facilement les r\'esultats 
du th\'eor\`eme \ref{class}, en interpr\'etant les r\'esultats de classification des $(\mathfrak{g},K)$-modules pour ${\rm Ind}_{MAN}^G(0\otimes \nu+\rho)$ 
avec l'\'egalit\'e $\nu+\rho=n\l$. Par exemple, l'assertion 
$\nu\in]-\frac{n}{2}, \frac{n}{2}[\setminus\{0\}$ \'equivaut \`a $\nu+\rho\in]0, n[\setminus\{\frac{n}{2}\}$, c'est \`a dire $\l\in]0, 1[\setminus\{\frac{1}{2}\}$. 
Les autres assertions d'unitarisabilit\'e s'obtiennent de la m\^eme mani\`ere.

\medskip

Le th\'eor\`eme \ref{class} est d\'emontr\'e.

\medskip

\begin{rema}
Le cas $n=1$ peut \^etre d\'eduit directement de la classification des repr\'esentations de $\SL(2,\bbR)$: adoptant la d\'emarche de M.~Lang dans
 \cite{Lan}, 
il convient \`a cet effet d'observer que  
l'espace des vecteurs $K$-finis de $\cF_\l(\mathbb{S}^1$) est la somme $\bigoplus H_{2l}$, $l\in\mathbb{Z}$, o\`u $H_m$ 
est l'espace de la repr\'esentation de $\SO(2)\subset\SL(2,\bbR)$ de caract\`ere 
$$
\chi_m:\left[\begin{array}{cc}\cos\theta&\sin\theta\\-\sin\theta&\cos\theta\end{array}
\right]\mapsto{\rm e}^{im\theta}.
$$
\end{rema}

\bigskip

{\bf Remerciements:} L'auteur tient \`a remercier vivement Thierry Levasseur, Alain Guichardet et Valentin Ovsienko pour 
leur grande disponibilit\'e, ainsi que Ranee Brylinski, Patrick Delorme et Pierre Lecomte pour d'enrichissantes discussions.
\backmatter


\end{document}